\documentclass[11pt]{article}
\usepackage{a4wide}

\usepackage{amsmath,amssymb,amsthm}
\usepackage{bbm}

\newtheorem{lemma}{Lemma}[section]
\newtheorem{proposition}[lemma]{Proposition}

\newtheorem{theorem}[lemma]{Theorem}
\newtheorem{corollary}[lemma]{Corollary}

\theoremstyle{definition}
\newtheorem{definition}[lemma]{Definition}

\newtheorem{example}[lemma]{Example}
\newtheorem{remark}[lemma]{Remark}

\usepackage[all]{xy}

\newcommand{\Real}{\mathbb{R}}
\newcommand{\Nat}{\mathbb{N}}

\newcommand{\Prob}{\mathsf{P}}
\newcommand{\Exp}{\mathsf{E}}

\newcommand{\Var}{\mathsf{Var}}

\newcommand{\Cov}{\mathsf{Cov}}
\newcommand{\COV}{\mathsf{COV}}

\newcommand{\RVeps}{\underline{\varepsilon}}
\newcommand{\RVdel}{\underline{\delta}}

\newcommand{\One}{\mathbbm{1}}

\newcommand{\Brown}{W}

\newcommand{\VRV}[1]{\mathbf{#1}}

\begin{document}

\title{Stochastic Calculus via Stopping Derivatives}

\author{Alex Simpson}

\maketitle

\begin{abstract}
We show that a substantial portion of stochastic calculus can be developed along similar lines to ordinary calculus, with derivative-based concepts driving the development. 
We define a notion of \emph{stopping derivative}, which is a form of right derivative with respect to stopping times. 
Using this, we define the \emph{drift} and \emph{variance rate} of a process as stopping derivatives for (generalised) conditional expectation and conditional variance respectively. Applying elementary, derivative-based methods, we derive a calculus of rules describing how drift and variance rate transform under operations on processes, culminating in a version of the multi-dimensional  Itô formula. Our approach connects with the standard machinery of stochastic calculus via a theorem 
establishing that continuous processes with zero drift coincide with random translations of continuous local martingales. 
This equivalence allows us to derive a  fundamental theorem of calculus for stopping derivatives, which relates the 
quantities of drift and variance rate, defined as stopping derivatives, 
 to parameters used in the description of a process as a stochastic integral.
\end{abstract}

\section{Introduction}
\label{section:intro}

In this article, we draw a new parallel between stochastic analysis and elementary real analysis, by developing an approach to the former along lines that closely follow the standard derivative-based development of the latter. In doing so, we address
a known conceptual issue with stochastic analysis as it is usually formulated: quantities that are  naturally derivative-like in conception, rather than being defined directly
in straightforward derivative-like terms, are   obtained only in circuitous ways.
We show that such quantities can be given simple, direct definitions using 
a new form of derivative, which we call a \emph{stopping derivative}. Stopping derivatives turn out to have a rich underlying theory that is
powerful enough to support  the development of significant parts of stochastic analysis 
along elementary derivative-based lines. The resulting 
development  not only parallels the 
standard development of real analysis, but in some ways directly generalises it. For example, we shall obtain a fundamental theorem of calculus in the stochastic setting (Theorem~\ref{theorem:stochastic-integral}) that includes a right-continuous version of the classical theorem as a special case.

To motivate the development that follows, let us consider, for illustrative purposes, the notion of \emph{drift} of a 
continuous-time real-valued stochastic process $(X_t)_t$, where by
 \emph{drift} we mean a measure of the rate of change of expected location. (Some authors call this \emph{drift rate} and instead use \emph{drift} to refer to the actual change in expected location.)
Let us  review two main ways in which this concept can be defined in conventional 
stochastic calculus. Our purpose in doing so is to underline the point that 
the definitions are indirect and technically involved. 
(Readers who do not have the technical background for the  discussion that follows
can take this point on trust and  skip the remainder of this paragraph.)
Firstly, in the case that $(X_t)_t$ is an It\^{o} process, its drift is  directly supplied as a parameter $(D_t)_t$  used in the description of $(X_t)_t$  in integral form,
\begin{equation}
\label{equation:integral-description}
X_t ~ = ~ X_0  + \int_0^t D_s\, ds + \int_0^t \sigma_s\,dW_s \enspace ,
\end{equation}
where the right-hand integral is an It\^{o} integral with respect to a given Brownian motion $(W_t)_t$.
If $(D_t)_t$ is right-continuous at time $s$, 
then the random variable $D_s$ is determined by the process $(X_t)_t$ alone, and $D_s$  can be said to be  the
\emph{drift of $(X_t)_t$ at time $s$}. 
The same notion can alternatively be defined, in a  more general context, by appealing to semimartingale theory. 
Any continuous semimartingale $(X_t)_t$
uniquely decomposes as
$(F_t + G_t)_t$, where $(G_t)_t$ is a continuous local martingale with initial value $0$ and $(F_t)_t$ is a (necessarily continuous) finite-variation process. 
In the case that $(F_t)_t$ has right derivatives along paths at time $s$, with these 
derivatives given by the random variable $D_s$,
then $D_s$ may again be considered to be the \emph{drift of $(X_t)_t$ at time $s$}. Importantly, in the special case of an 
 It\^{o} process, the second approach to defining drift agrees with the first. 

In contrast with  the technical definitions above,
there is a natural intuitive notion of drift that is manifestly simple in conception.
The \emph{drift of a process $(X_t)_t$} at time $s$ is intended to  quantify
the \emph{rate of change}, at time $s$, \emph{of the expected position} of the process, given that its actual position at time $s$ is known.
This idea is often encountered in an informal motivational role (as in, e.g..~\cite[\S5.1]{RWII}), but it is not usually taken as a template for a mathematical definition. (One exception to this,
which we discuss further in Section~\ref{section:end-notes}, is 
Nelson \cite[Chapter 11]{Nelson}.)
In fact several difficulties arise if one  does try to turn this naive notion of drift into a precise definition using
the most obvious mathematical translation; that is, if one defines \emph{the
drift at time $s$} as a right derivative at $t = s$ for the \emph{conditional expectation} function
\begin{equation}
\label{rd:exp}
 t ~ \mapsto ~ \Exp[X_t \mid \mathcal{F}_s] \enspace ,
 \end{equation}
where we assume that $(X_t)_t$ is adapted to a given filtration $(\mathcal{F}_t)_t$.
(One reason for asking for a \emph{right} derivative is that the conditioning on $\mathcal{F}_s$ makes only the future behaviour of $(X_t)_t$ uncertain.)
Interpreting this literally, one obtains a definition of drift at time $s$ as
\begin{equation}
\label{equation:drift-as-rd}
D_s ~ := ~ \lim_{t \downarrow s} \, \frac{\Exp[X_t \mid \mathcal{F}_s] - \Exp[X_s \mid \mathcal{F}_s]}{t-s} \enspace ,
\end{equation}
where the right-hand side expresses the  limit as $t$ tends to $s$ from above. 
(Of course,~\eqref{equation:drift-as-rd} could be written more simply, because
$\Exp[X_s \mid \mathcal{F}_s] = X_s$, but we leave it in unsimplified form to emphasise the general derivative formula involved.)
In order to make sense of~\eqref{equation:drift-as-rd},
we
need to agree on what kind of limit of random variables we are talking about (convergence in probability, almost sure convergence, \dots). 
There is no clear  canonical choice, but the discussion that follows  applies irrespective of the convergence notion chosen.

There are good reasons that one does not typically encounter~\eqref{equation:drift-as-rd} as a definition of drift in the literature. Firstly, one can make sense 
of~\eqref{equation:drift-as-rd} only if $\Exp[X_t \mid \mathcal{F}_s]$ is defined for $t$ sufficiently close to $s$. Using the standard definition of conditional expectation, this requires 
$X_t$, for such $t$,  to be integrable, which loses
the generality of It\^o processes. If one instead makes use of the
\emph{generalised conditional expectation} of Meyer~\cite{meyer}, 
then $X_t$ does not have to be integrable for $\Exp[X_t \mid \mathcal{F}_s]$ to be defined. 
However,
the generalised conditional expectation 
$\Exp[X_t \mid \mathcal{F}_s]$  itself need not be defined for It\^o processes,
so some loss of generality 
remains.
One might, in principle, be prepared to accept such loss of generality, if definition~\eqref{equation:drift-as-rd} turned out to have other advantages. However, it does not behave well at all. 
Most damningly, the naive drift defined by~\eqref{equation:drift-as-rd} 
has very little to do with the established concept of drift, as reviewed above for semimartingales (including It\^o processes).
%
%
For example, either one of naive drift or semimartingale drift may be defined and the other not. Moreover, even if both are defined, they need not agree with each other. 
In the face of such disagreement, the semimartingale definition has the distinct advantage of fitting into the well-developed framework of stochastic calculus, with its useful toolkit of rules for manipulating stochastic processes. In contrast, definition~\eqref{equation:drift-as-rd} does not seem compatible  with any such practical calculus.  

The key observation underpinning the present article is that, by making one simple and conceptually natural modification
to~\eqref{equation:drift-as-rd} above, one does obtain a  definition of drift that both  directly formalises
the idea that drift is rate-of-change of expected position, and which simultaneously generalises the
semimartingale definition. 
The modification we make to~\eqref{equation:drift-as-rd} is that we \emph{randomise time}. That is, we let the time variable $t$ range over \emph{stopping times} rather than deterministic times. We thus define the \emph{drift 
at time $s$} as a right derivative  at $T = s$
for the function
\begin{align}
\label{equation:E-T-s}
 & T ~ \mapsto ~ \Exp[X_T \mid \mathcal{F}_s] 
\end{align}
where $T$ ranges over finite stopping times with $T \geq s$.
Formally, this means defining
\begin{equation}
\label{equation:drift-as-sd}
D_s ~ := ~ \lim_{T \downarrow s} \, \frac{\Exp[X_T \mid \mathcal{F}_s] - \Exp[X_s \mid \mathcal{F}_s]}{\,\Exp[T-s \mid \mathcal{F}_s]\,} \enspace ,
\end{equation}
where the limit is over finite stopping times $T > s$, and the  conditional expectation
in the denominator ensures the $\mathcal{F}_s$-measurability of the quotient.
We call such right derivatives, defined with respect to stopping times,  \emph{stopping derivatives}.


All the problems besetting the original deterministic-time definition~\eqref{equation:drift-as-rd}  are resolved by the move to a stopping derivative~\eqref{equation:drift-as-sd}. First, 
as Theorem~\ref{theorem:stopping-limit}, we shall establish that there is a canonical  convergence notion
for limits of the form~\eqref{equation:drift-as-sd} above.
There is accordingly no debate over what form of convergence  to use.
Next, the problems concerning the potential undefinedness of  conditional expectation do not arise.
Under a mild path-continuity assumption, 
 $\Exp[X_T \mid \mathcal{F}_s]$ 
is always defined for sufficiently small $T$, as long as conditional expectation is interpreted 
in the generalised sense of Meyer~\cite{meyer}. This follows from Propositions~\ref{proposition:exp-continuous} 
below,  which proves 
that the function~\eqref{equation:E-T-s} 
satisfies a stronger property called \emph{stopping continuity}.
Finally, in the case of processes that are continuous semimartingales,
the stopping-derivative notion of drift
does coincide with the established notion. 
Moreover, since the definition of drift 
as a stopping derivative is not 
dependent on any form of process decomposition, 
it is potentially applicable more widely to processes that are not  semimartingales. 

An important feature of stopping derivatives is that they
give rise to their own calculus of rules governing how they  transform under operations on processes. 
Such rules offer an elementary differential-calculus-based approach to stochastic calculus. 
To give a flavour of this, let us define a second fundamental quantity as a stopping derivative. 
Given that the canonical measures of the location and spread of a probability distribution are
expectation and variance, and we have defined drift as a stopping derivative for conditional expectation, it is natural to 
consider, as a companion quantity, a stopping derivative for the \emph{conditional variance} function
\begin{equation}
\label{rd:var}
  T ~ \mapsto ~ \Var[X_T \mid \mathcal{F}_s] := \Exp[(X_T - \Exp[X_T \mid \mathcal{F}_s])^2 \mid \mathcal{F}_s] \enspace .
\end{equation}
The resulting quantity, which we call \emph{variance rate}, is then
defined by simply replacing conditional expectation with conditional variance in the numerator of  \eqref{equation:drift-as-sd}, \emph{viz}
\begin{equation*}
\label{equation:varrate-as-sd}
\lim_{T \downarrow s} \, \frac{\Var[X_T \mid \mathcal{F}_s]}{\,\Exp[T-s \mid \mathcal{F}_s]\,} \enspace ,
\end{equation*}
where we have simplified the numerator using the equality  $\Var[X_s \mid \mathcal{F}_s] = 0$.

Having the two quantities of drift and variance rate in combination, one can obtain  derivative-based formulations of key results from stochastic calculus. 
For example, as Theorem~\ref{theorem:ito-formula},
we shall derive It\^o's   formula  in (a generalisation of) the following form. 
Suppose that $(X_t)_t$  has drift $D_s$ and 
variance rate $A_s$ at time $s$,  and suppose that 
$f$ is twice continuously differentiable. Then the process $(f(X_t))_{t}$ has drift $E_s$ and variance rate $C_s$, at time $s$, where:
\begin{equation}
\label{equation:ito}
E_s = f'(X_s)D_s + \frac{f''(X_s)}{2}A_s  \qquad \text{and}  \qquad C_s = (f'(X_s))^2A_s \enspace .
\end{equation}
The development  of stochastic calculus, based on drift and variance rate,  turns out to be remarkably elementary, with the rules of calculus falling out straightforwardly
from the formal definitions as  stopping derivatives. The entire development is carried out using the epsilon-delta style of reasoning familiar from basic real analysis, albeit with some unavoidable added complexity arising from the
need to manipulate stopping times and to respect measurability constraints. 
Stochastic integrals play no part in the development.

One further  important feature of the  approach outlined above is that our definitions of
drift $D_s$ and variance rate $A_s$ define these values \emph{locally} at  an individual  time $s$. 
In fact, 
the definitions enjoy a stronger localisability property: they easily adapt to define
the drift $D_S$ and variance rate $A_S$  locally at any  finite stopping time $S$. 
It is then natural to say that a process has \emph{zero drift} if it  has drift $0$ at all finite stopping times. 
As Theorem~\ref{theorem:driftless-martingale}, we prove that 
a continuous process  has zero drift  if and only if it is a random translation of a local martingale. That is, for continuous processes, the stopping-derivative notion of driftlessness is intimately  connected to  the established martingale-based one. 
This key equivalence forms the basis for establishing, as  Theorem~\ref{theorem:stochastic-integral}, a fundamental theorem of calculus for stopping derivatives.  Suppose 
$(X_t)_t$ is defined via an It\^o-style  stochastic integral~\eqref{equation:integral-description}, using a Brownian motion $(W_t)_t$ as the integrator,
and suppose that the integrand processes $(D_s)_s$ and $(\sigma_s)_s$ are right continuous. Then
$(D_s)_s$ and $({\sigma_s}^2)_s$ can be recovered respectively  as the stopping-derivative-defined drift process and variance-rate process for $(X_s)_s$. Thus, for  It\^o processes, 
our stopping-derivative-based definitions provide an inverse to stochastic integration, where the derivative-based and integral-based notions of drift coincide, and the derivative-based notion of variance rate recovers the square of the diffusion coefficient $\sigma_s$ in the integral.

The article is structured as follows. 
After a brief review of background material in Section~\ref{section:preliminaries},
Section~\ref{section:stopping-functionals}  introduces technical machinery needed to define stopping derivatives: 
specifically  \emph{stopping functionals} (of which conditional expectation and conditional variance are examples) and their canonical 
\emph{stopping limits}. 
The notion of \emph{stopping derivative} itself is defined in Section~\ref{section:stopping-derivatives} as a particular 
stopping limit, and  its basic properties are  established.
Sections~\ref{section:drift} and~\ref{section:diffusion} then examine the concepts of 
\emph{drift} and \emph{variance rate}, defined as stopping derivatives for
conditional expectation and conditional variance respectively, ending with 
our version of the It\^o formula. Section~\ref{section:covariance} considers  \emph{covariance rate} for pairs of processes, leading to the multi-dimensional  It\^o formula.
Section~\ref{section:fundamental} addresses the topic of local martingales from the viewpoint of stopping derivatives.
In particular, we 
prove that a continuous process has zero drift if and only if it is a random translation of a local martingale.
In Section~\ref{section:integration-uniqueness},
we prove a  fundamental theorem of calculus for stopping derivatives. 
We also prove a uniqueness result establishing sufficient conditions under which 
a continuous zero-drift process is  determined by its variance rate.
Finally, in Section~\ref{section:end-notes}, we set our work in perspective with a discussion of some of the 
other derivative-based approaches to quantities similar to 
drift and variance rate that one finds in the literature. We also discuss possible directions for further investigation.

\section{Preliminaries}
\label{section:preliminaries}

We consider continuous-time real-valued processes $(X_t)_{t \geq 0}$ adapted to a filtration $(\mathcal{F}_t)_{t \geq 0}$ over a sample space $\Omega$ with $\sigma$-algebra $\mathcal{E}$ and probability measure $\Prob$.
We often write a process $(X_t)_{t \geq 0}$  as just $(X)$. When we mention a process $(X)$ and filtration $(\mathcal{F})$ in the same context, the process is assumed to be adapted, unless explicitly stated otherwise. Further assumptions are imposed explicitly as and when needed. In particular, we do not make any additional default assumptions on the filtration; so we do not, for example, assume that it satisfies the \emph{usual conditions}.
In the measurability requirements on processes, measurability 
over $\Real$ is with respect to the Borel $\sigma$-algebra $\mathcal{B}$. If $\mathcal{G} \subseteq \mathcal{E}$ is a $\sigma$-subalgebra, then by a \emph{$\mathcal{G}$-partition of $\Omega$} we mean 
a family $\{P_i \mid i \in I\} \subseteq \mathcal{G}$ satisfying $\bigcup_{i \in I} P_i = \Omega$ and 
$P_i \cap P_j = \emptyset$ whenever $i \neq j$. 

We shall often need  processes $(X)$  to be  \emph{progressive}, allowing us to evaluate $X$ at a finite stopping time $T$ to obtain an $\mathcal{F}_T$-measurable $X_T$. (In this work, we use capital $S, T, U,\dots$ for stopping times, and lower case $s,t,u,\dots$ for deterministic times.)
In the majority of  cases, we obtain progressiveness indirectly via path-continuity  properties. We say that an (adapted by default) process 
$(X)$ is \emph{right continuous}
or \emph{cadlag} or  \emph{left continuous} or \emph{continuous} if, for all $\omega \in \Omega$, the path $t \mapsto X_t(\omega)$  enjoys the stated property. As is well known, right-continuous processes are automatically \emph{optional} (hence progressive) and left-continuous processes are even \emph{previsible} (hence optional). Further, any cadlag process $(X)$ determines a derived left-continuous process $\omega,t \mapsto X_{t-}(\omega)$, where
\[
X_{t-}(\omega) ~ : = ~ \lim_{s \uparrow t} X_s(\omega) \enspace ,
\]
i.e., the limit of $X_s(\omega)$ as $s$ tends to $t$ from below. We frequently exploit this notation in contexts in which we have a  finite stopping time $T$ bounded below by some distinguished stopping time $S \leq T$. 
In such a situation, we write $X_{T-}$ to mean
\begin{equation}
\label{eqn:T-minus}
X_{T-}(\omega) ~ := ~ 
 \begin{cases} 
 X_{T(\omega)-} & \text{if $T(\omega) > S(\omega)$} \\
 X_{T(\omega)} & \text{if $T(\omega) = S(\omega)\,$.}
 \end{cases}
\end{equation}

We also frequently exploit the cadlag property to define stopping times. For example, if $(X)$ is a cadlag process, $S$ a finite stopping time and $D_S$  a positive $\mathcal{F}_S$-measurable random variable, then the \emph{first approach} beyond $S$ to the random set
$\{x \in \Real \mid |x - X_S| \geq D_S\}$ is the stopping time: 
\[
T(\omega) ~ := ~ \inf \{  t > S(\omega) \mid \, \text{$|X_t(\omega) - X_{S}(\omega)| \geq D_S(\omega)$ or $|X_{t-}(\omega) - X_{S}(\omega)| \geq D_S(\omega)$} \} \enspace .
\]
In similar definitions in the sequel, we shall abbreviate the right-hand side using $\pm$ in the notation, \emph{viz}:
\begin{equation}
\label{equation:plus-minus}
\inf \{  t > S(\omega) \mid \, |X_{t\pm}(\omega) - X_{S}(\omega)| \geq D_S(\omega) \} \enspace .
\end{equation}

Path-differentiability properties will also play a role, in particular the property of \emph{right differentiability}. We shall need the right-differentiable formulation of the fundamental theorem of calculus in the form: if a measurable function $f$ on $[a,b)$ satisfies $\int_a^b |f(s)| {\mathit{d}s} < \infty$ and if $f$ is right continuous at $c \in [a,b)$ then the continuous function $t \mapsto \int_a^t f(s) {\mathit{d}s}$ is right differentiable at $c$ with right derivative $f(c)$.

Finally, in the sequel, we shall often assert (in)equalities between random expressions defined  up to almost-sure equality. In such cases, unless something different is specified, we mean that the relation in question holds almost surely.


\section{Stopping functionals and stopping limits}
\label{section:stopping-functionals}

As adumbrated in the introduction, we are interested in  an appropriate form of right derivative, with respect to  stopping times $T$, for the operations that return \emph{conditional expectation} and \emph{conditional variance} for a process $(X)$. We also wish  to define such derivatives locally at finite stopping times $S$. Accordingly, we require $(X)$ to be progressive and we consider the functions that map 
pairs of finite stopping times $S \leq T$ to
\begin{align}
\label{def:exp}
& \Exp[X_T \mid \mathcal{F}_S] \\
\intertext{and} 
\label{def:var}
& \Var[X_T \mid \mathcal{F}_S]  ~ := ~ \Exp[(X_T - \Exp[X_T \mid \mathcal{F}_S])^2  \mid \mathcal{F}_S] \enspace .
\end{align}
These are, in general, partial functions, since the above expressions need not be defined for all pairs $S \leq T$ of finite stopping times.

The functions~\eqref{def:exp} and~\eqref{def:var} above are examples of what we call \emph{stopping functionals}, the general definition of which requires only 
a filtration $(\mathcal{F}_t)_{t \in [0,\infty)}$ on a sample space $\Omega$. 


\begin{definition}[Stopping functional]
A   \emph{stopping functional} is a partial function $F$ that takes as argument a pair of stopping times  $S \leq T < \infty$, and which returns, if defined, an $\mathcal{F}_S$-measurable $\Real$-valued random variable 
 determined up to $\mathcal{F}_S$-almost-sure equality. 
 We write $F[T \mid \mathcal{F}_S]$  for the application of $F$ to the pair $(S,T)$, since this notation is suggestive both of the measurability requirement and also of $S$ playing a conditioning role consistent with the  \emph{partition property} below, which we further require any stopping functional $F$ to satisfy. 
 \begin{quotation}
 \noindent
 \emph{Partition property:} 
 for any stopping time $S < \infty$, countable $\mathcal{F}_S$-partition $\{P_i \mid i \in I\}$ of $\Omega$,
 and family $(T_i)_{i \in I}$ of stopping times with $S \leq T_i < \infty$ for all $i$, if
 both the right-hand and left-hand expressions below are defined then they are equal:
\begin{equation}
\label{equation:locality}
F\left[\left. \sum_{i \in I} \, \One_{P_i } \cdot T_i \, \right| \, \mathcal{F}_S \right] ~ = ~  \sum_{i \in I} \, \One_{P_i } \cdot  F[T_i \mid \mathcal{F}_S] \enspace .
\end{equation}
\end{quotation}
\noindent
We shall also consider the following strengthening of the partition property, which we do not assume as a default property of stopping functionals. 
\begin{quotation}
\noindent
\emph{Strong partition property:} If the right-hand expression of~\eqref{equation:locality} is defined then so is the 
left-hand expression and they are equal.
 \end{quotation}
 \end{definition}
 
 It is easy to check that conditional expectation~\eqref{def:exp} and conditional variance~\eqref{def:var} both satisfy the partition property and are therefore examples of stopping functionals. If we adhere to the standard notion of conditional expectation $\Exp[Y \mid \mathcal{G}]$, in which $Y$ is required to be integrable, then~\eqref{def:exp} and~\eqref{def:var} do not satisfy the strong partition property. However, if we instead adopt the generalised form of conditional independence introduced by Meyer~\cite{meyer} and studied systematically in~\cite[\S{I.4}]{HWY}, then  the strong partition property does hold. Recall that, in its generalised form, $\Exp[Y \mid \mathcal{G}\,]$ is defined 
 precisely when $\Exp[\,|Y| \mid \mathcal{G}\,]$ is almost surely finite (where $\Exp[\,|Y| \mid \mathcal{G}\,]$ is the 
 $[0,\infty]$-valued conditional expectation defined for arbitrary non-negative random 
 variables), in which case one defines
 \[
 \Exp[Y \mid \mathcal{G}] ~ := ~ \Exp[Y^+ \mid \mathcal{G}] - \Exp[Y^- \mid \mathcal{G}] \enspace ,
 \]
 using two non-negative conditional expectations on the right. 
 This definition, with its separation of positive and negative parts, is manifestly natural, being as it is 
a direct conditional analogue of a standard definition of ordinary unconditional expectation.  
 
 One of the properties of generalised
conditional expectation we shall make use of  is that the identity $\Exp[Y \mid \mathcal{G}] = Y$ holds for arbitrary $\mathcal{G}$-measurable $Y$, without any integrability conditions. 
More generally, for any $\mathcal{G}$-measurable random variables $Y, \RVdel$ with $\RVdel \geq 0$, and any 
random variable $Z$ such that $|Z - Y| \leq \RVdel$, it holds that $\Exp[Z \mid \mathcal{G}]$ is defined
and satisfies $|\Exp[Z \mid \mathcal{G}] - Y| \leq \RVdel$.
(These properties are first used in the proof of Proposition~\ref{proposition:exp-continuous} below.) 
Another property that holds for generalised conditional expectation is that, for any
countable $\mathcal{G}$-partition $\{P_i \mid i \in I\}$ of $\Omega$, and family of random variables
$(Y_i)_{i \in I}$, if the generalised conditional expectation $\Exp[Y_i \mid \mathcal{G}]$ is defined, for all $i$, then
\begin{equation}
\label{equation:countable-sum}
\Exp\left[\left. \sum_{i \in I} \, \One_{P_i } \cdot Y_i \, \right| \, \mathcal{G} \right] ~ = ~  \sum_{i \in I} \, \One_{P_i } \cdot  \Exp[Y_i \mid \mathcal{G}] \enspace ,
\end{equation}
in particular, the left-hand side is defined. 
(This will be important in the proof of Theorem~\ref{theorem:martingale-driftless} below.)
This implies that the strong partition property 
holds for the generalised-conditional-expectation stopping functional~\eqref{def:exp};
and hence also for conditional variance~\eqref{def:var}, as long as 
generalised conditional expectation is used in the definition of the latter.
Since properties such as the aforegoing  are needed throughout the article,
 we henceforth adopt generalised conditional expectation  as our default notion 
 of {conditional expectation}, and we shall use the standard conditional-expectation  notation
 $\Exp[Y \mid \mathcal{G}]$ for it without further explanation. 
 

Returning to the general theory of stopping functionals, limits of the form below, for the conditional expectation stopping functionals $F$ and finite stopping times $S$, will play a fundamental role in the sequel,
 \begin{equation}
 \label{equation:stopping-limit} 
 \lim_{T \downarrow S}\,  F[T \mid \mathcal{F}_S] \enspace .
 \end{equation}
 The notation means the limit of $F[T \mid \mathcal{F}_S]$ as the stopping time $T$ tends to $S$ from above. This is a limit
 of a net of $\mathcal{F}_S$-measurable random variables indexed by the set
 \begin{equation}
 \label{equation:stopping-net}
 \{\text{$T$ a finite stopping time} \mid T > S\} \enspace ,
 \end{equation}
which is filtered under the operation of pairwise minimum. The required notion of limit  is  
therefore a form of net convergence for $\mathcal{F}_S$-measurable random variables. 
The definition below,
considers several possible definitions for such forms of limit  $L_S = \lim_{T \downarrow S}\, F[T \mid \mathcal{F}_S]$.

\begin{definition}[Forms of net convergence] 
\label{definition:net-convergence}
Suppose $F$ is a stopping functional, $S$ a finite stopping time and  $L_S$ an $\mathcal{F}_S$-measurable
 random variable.

\begin{enumerate}

\item 
\label{s-l:c}
\emph{(Convergence in probability)} 
We say that $F[T \mid \mathcal{F}_S]$ converges \emph{in probability} to $L_S$ if,
for 
any  real $\varepsilon >0$, it holds that
\[
\lim_{T \downarrow S} \, \Prob[| F[T \mid \mathcal{F}_S] - L_S | \leq  \varepsilon] ~ = ~ 1 \enspace ,
\]
where the limit is convergence in $\Real$ over the net of stopping times $\{T < \infty \mid T > S\}$. In expanded form, this means that,  
for every $\delta > 0$, there exists
 some finite stopping time $U > S$ such that, for all stopping times $T$ with $S < T \leq U$, it holds that
 \begin{equation}
 \label{equation:s-l:c}
 \Prob[| F[T \mid \mathcal{F}_S] - L_S | \leq  \varepsilon] ~ \geq ~ 1- \delta  \enspace .
 \end{equation}

\item 
 \label{s-l:d}
 \emph{(Persistent convergence in probability)}
 We say that $F[T \mid \mathcal{F}_S]$  converges \emph{persistently in probability} to $L_S$ if,
 for any  reals $\varepsilon, \delta  >0$, there exist $W \in \mathcal{F_S}$ and a finite stopping time
 $U > S$ such that, for all stopping times $T$ with $S < T \leq U$, it holds that
 $\Prob(\{\omega \in W \mid \, | F[T \mid \mathcal{F}_S](\omega) - L_S(\omega) | \leq  \varepsilon\}) = \Prob(W) \geq 1 - \delta$.

 \item 
 \label{s-l:b}
 \emph{(Uniform convergence)} We say that $F[T \mid \mathcal{F}_S]$ converges \emph{uniformly} to $L_S$ if,
 for 
 any  real $\varepsilon >0$, there exists
 some finite stopping time $U > S$ such that, for all stopping times $T$ for which $S < T \leq U$, it holds that
 \[
 | F[T \mid \mathcal{F}_S] - L_S | ~ \leq ~ \varepsilon \enspace .
 \]

\item
 \label{s-l:a}
 \emph{(Randomly uniform convergence)} We say that $F[T \mid \mathcal{F}_S]$ converges \emph{randomly uniformly} to $L_S$ if,
 for any $\mathcal{F}_S$-measurable random variable $\RVeps>0$, there exists
 some finite stopping time $U > S$ such that, for all stopping times $T$ for which  $S < T \leq U$, it holds that
 \[
 | F[T \mid \mathcal{F}_S] - L_S | ~ \leq ~ \RVeps \enspace .
 \]

   \end{enumerate}
 \noindent
  In all four properties above, there is a tacit assumption that $F[T \mid \mathcal{F}_S]$ is defined for all
  stopping times $T$ for which $S < T \leq U$. It need not be defined for other stopping times.

\end{definition}

The four defined notions of convergence all straightforwardly generalise to the more general context of
convergences $\lim_{i \in I} X_i$ for families $(X_i)_{i \in I}$ of random variables indexed by filtered partial orders $I$. We do not present them in this more general form, because  we shall only be concerned with convergences of the form~\eqref{equation:stopping-limit}.
Nevertheless, some remarks about the general notions are in order. The notion of 
\emph{convergence in probability} is standard, and takes the expected form. 
The property we call \emph{persistent convergence in probability} is included as a replacement for \emph{almost-sure convergence}.
It is equivalent to almost-sure convergence in the special case of countable index sets $I$ (so, in particular, in the case of sequence convergence). In the case of uncountable index sets, however, such as our index set~\eqref{equation:stopping-net}, the notion of almost sure convergence is inappropriate because it is not invariant under choice of versions of random variables defined only up to almost-sure equality.
 The property of \emph{uniform convergence} takes the expected form. The property we  call  \emph{randomly uniform convergence} is a strengthening of uniform convergence that is interesting only in the case of sufficiently rich index sets $I$, for  it can  be satisfied only in trivial ways when $I$ is countable.  For our index set of stopping times~\eqref{equation:stopping-net}, this strong notion of convergence turns out to be particularly fruitful. 

%
%
%

It is an easy observation that the implications  \ref{s-l:a} $\implies$ \ref{s-l:b} $\implies$ \ref{s-l:d} $\implies$ \ref{s-l:c} hold between the notions of convergence introduced in Definition~\ref{definition:net-convergence}. Moreover, these implications hold 
in the general context of convergence of nets of random variables $\lim_{i \in I} X_i$. It is also easy to see that the converse implications do not hold  for general nets of random variables. Somewhat surprisingly, for the limits we are interested in, i.e., those of the form~\eqref{equation:stopping-limit} above, the implications  \emph{do}  reverse.
 There is therefore just one canonical notion  of convergence for stopping-time-indexed nets of random variables $F[T \mid \mathcal{F}_S]$ given by
 a stopping functional $F$.

 \begin{theorem}
 \label{theorem:stopping-limit}
 Let $F$ be a stopping functional, $S$ a finite stopping time and $L_S$ an $\mathcal{F}_S$-measurable random variable
 Then the four notions of convergence
 in Definition~\ref{definition:net-convergence} coincide.  \end{theorem}

 \begin{proof}
  As remarked above, the implications \ref{s-l:a} $\implies$ \ref{s-l:b} $\implies$ \ref{s-l:d} $\implies$ \ref{s-l:c} are obvious. We complete the proof by establishing the chain of reverse implications.
 
 To prove \ref{s-l:c} $\implies$ \ref{s-l:d}, assume \ref{s-l:c} and let  $\varepsilon, \delta  >0$ be arbitrary.
 Let $U$ be such that, for all  $T < \infty$ such that  $S < T \leq U$, we have that $F[T \mid \mathcal{F}_S]$ is defined 
 and~\eqref{equation:s-l:c} holds . For any
 $T < \infty$ with $S < T \leq U$, define
 \begin{equation}
 \label{equation:def-VT}
 V_T ~ := ~ \{\omega \in \Omega \mid \, |F[T \mid \mathcal{F}_S](\omega) - L_S(\omega)| > \epsilon\} \enspace .
 \end{equation}
 Clearly $V_T \in \mathcal{F}_S$ and $\Prob(V_T) < \delta$. 
 The \emph{essential union} of all $V_T$, is a set $V \in \mathcal{F}_S$ of smallest probability such that, for every $T$, we have 
 $\Prob(V_T \backslash V) = 0$.
 By a standard construction, we can obtain this as a countable disjoint union
  $V := \bigcup_n V_n$, where each $V_n$ is taken from 
 the closure of the set 
 $\{V_T \mid \text{$T < \infty$ and $S < T \leq U$}\}$ under finite intersections and relative complements.
 For every $n$, let $T_n$ be a finite stopping time such that $S < T_n \leq U$ and $V_n \subseteq V_{T_n}$. Define the stopping time $U'$ by:
 \[ U'(\omega) ~ := ~ \begin{cases} T_n(\omega) & \text{if $\omega \in V_n$} \\
 U(\omega) & \text{if $\omega \notin V$} \enspace .
 \end{cases}
 \]
 Clearly $U' \leq U$. 
 
 Suppose $T$ is a finite stopping time with $S < T \leq U'$. 
Because  $\Prob({V_T} \backslash V) =0$, it follows from~\eqref{equation:def-VT} that
\begin{equation}
\label{equation:show-three}
\Prob(\{\omega \in V^C \mid \, |F[T \mid \mathcal{F}_S](\omega) - L_S(\omega)| \leq \epsilon\}) ~ = ~ \Prob(V^C) \enspace .
\end{equation}
Define $W := V^C$. Below we show that $\Prob(W) \geq 1 - \delta$. Hence, by~\eqref{equation:show-three} above,
$W$ and $U'$ satisfy the properties needed to show 
 persistent convergence in probability.
 
 It remains to show that $\Prob(W) \geq 1 - \delta$. By the definition of $U'$, we have
 \begin{align}
 \nonumber
 F[ U' \mid \mathcal{F}_S] ~  & = ~ F \left[ \left. \One_{W} \cdot U + \sum_{n} \One_{V_n} \cdot T_n\,  \right| \, \mathcal{F}_S\right]
 \\
 \label{equation:WVn-sum}
 & = ~ \One_{W} \cdot F[ U \mid \mathcal{F}_S] + \sum_{n} \One_{V_n} \cdot F[T_n \mid  \mathcal{F}_S]
 & & \text{by the partition property},
 \end{align}
 where the partition property applies because both sides are defined since all stopping times are $\leq U$.
 We can now reason as follows.
 \begin{align*}
 1 - \delta ~ & \leq ~ \Prob[ | F[ U' \mid \mathcal{F}_S] - L_S| \leq \varepsilon]  & & \text{by~\eqref{equation:s-l:c}} \\
& = ~  \Prob(\{ \omega \in W  \mid \, |F[ U \mid \mathcal{F}_S](\omega)  - L_S(\omega)| \leq \varepsilon\}) \\
 & ~~~ ~~~~~~ + \sum_n \Prob(\{\omega \in V_n \mid \, |F[ T_n \mid \mathcal{F}_S](\omega)  - L_S(\omega)| \leq \varepsilon\})
 & & \text{by~\eqref{equation:WVn-sum}} \\
& = ~ \Prob(\{ \omega \in W \mid \, |F[ U \mid \mathcal{F}_S](\omega)  - L_S(\omega)| \leq \varepsilon\})
\\ & \leq ~ \Prob(W)  \enspace ,
 \end{align*}
 where the second equality holds because it follows from~\eqref{equation:def-VT} and the inclusion $V_n \subseteq V_{T_n}$ that
 \[
 \Prob(\{\omega \in V_n  \mid \, |F[ T_n \mid \mathcal{F}_S](\omega)  - L_S(\omega)| \leq \varepsilon\}) ~ = ~ 0 \enspace .
 \]
 
 To prove \ref{s-l:d} $\implies$ \ref{s-l:b}, assume \ref{s-l:d} and let $\varepsilon > 0$ be arbitrary. Using property \ref{s-l:d}, for any $n > 0$, let $W_n \in \mathcal{F}_S$ and $U_n$ be such that $\Prob(W_n) \geq \frac{n}{n+1}$ and  $U_n > S$ is a finite stopping time satisfying, 
 for all stopping times $T_n$ with $S < T_n \leq U_n$,
 \begin{equation}
 \label{equation:from-sld}
 \Prob( \{\omega \in W_n \mid \, | F[T_n \mid \mathcal{F}_S](\omega) - L_S(\omega) | \leq  \varepsilon\}) ~ = ~ \Prob(W_n)
 \enspace .
 \end{equation}
 Obviously, we can assume $U_n \leq U_1$.
 For $n > 0$, define $V_n := \bigcup_{i = 1}^n W_i \backslash \bigcup_{i = 1}^{n-1} V_i$.
 The sets $(V_i)_{i \geq 1}$ are clearly pairwise disjoint and $W_n \subseteq \bigcup_{i = 1}^n V_i$.
 It follows that $\Prob\left(\bigcup_{i = 1}^n V_i\right) \geq \frac{n}{n+1}$, hence
 $\Prob(V) = 1$, where  $V :=  \bigcup_{i = 1}^\infty V_i$. Define a stopping time $U$ by:
 \[
 U(\omega) ~ := ~ \begin{cases} U_n(\omega) & \text{if $\omega \in V_n$} \\
 U_1(\omega)  & \text{if $\omega \not\in V$} \enspace .
 \end{cases}
 \]
 Clearly $S < U \leq U_1$. We show that $U$ satisfies the property required to establish   uniform convergence \ref{s-l:b}. Suppose that $T$ is a stopping time with $S < T \leq U$. We need to show that $ \Prob[\, | F[T \mid \mathcal{F}_S] - L_S | \leq  \varepsilon] = 1$.
 
 Define
 \[
 T_n(\omega) ~ := ~ \begin{cases}
   T(\omega) & \text{if $\omega \in V_n$} \\
    U_n(\omega) & \text{if $\omega \not\in V_n$}  \enspace .
   \end{cases}
 \]
 Because $T \leq U$, we have $T_n \leq U_n$. Hence by~\eqref{equation:from-sld} and the fact that $V_n \subseteq W_n$, it holds that
 \begin{equation}
 \label{equation:Vn}
 \Prob(\{\omega \in V_n \mid \, | F[T_n \mid \mathcal{F}_S](\omega) - L_S(\omega) | \leq  \varepsilon\}) = \Prob(V_n) \enspace .
 \end{equation}
 Because $\{V_n\}_n \cup \{V^C\}$ is an $\mathcal{F}_S$-partition of $\Omega$, we have, by the definition of the stopping times $T_n$,
  \begin{align*}
  & \Prob[ | F[T \mid \mathcal{F}_S] - L_S | \leq  \varepsilon]  \\
  & ~~~ = ~  \Prob\left[ \, \left|\,  F\left[ \left. \One_{V^C} \cdot T + \sum_n \One_{V_n} \cdot T_n \, \right| \,  \mathcal{F}_S\right] - L_S \right| \leq  \varepsilon\right]  \\[1ex]
  & ~~~ = ~ \Prob\left[ \, \left|\, \left(\One_{V^C} \cdot F[T \mid \mathcal{F}_S] + \sum_n \One_{V_n} \cdot F[T_n \mid \mathcal{F}_S] \right)  - L_S \right| \leq  \varepsilon \right] & & \text{by the partition property}
  \\[1ex]
  & ~~~ = ~ \sum_n\, \Prob(\{\omega \in V_n \mid \, |F[T_n \mid \mathcal{F}_S](\omega) - L_S(\omega)| \leq \varepsilon \})
  & & \text{because $\Prob(V^C) = 0$} \\
  & ~~~ = ~ \sum_n\, \Prob(V_n) & & \text{by~\eqref{equation:Vn}}
  \\
  & ~~~ = ~ P(V) ~ = ~ 1 & & \text{because $V = \biguplus_{i = 1}^\infty V_i$ \enspace .}
  \end{align*}
  Note that both sides are defined in the application of the partition property, because all the stopping times involved are $\leq U_1$.
  
 To prove \ref{s-l:b} $\implies$ \ref{s-l:a}, assume \ref{s-l:b} and let $\RVeps > 0$ be an $\mathcal{F}_S$-measurable random variable. For any $n \geq 1$, let $U_n$ be a finite stopping time, as given by uniform convergence \ref{s-l:b}, such that $U_n > S$ and,
 for all $T_n$ with $S < T_n \leq U_n$, it holds that  
 \begin{equation}
 \label{equation:one-on-n}
 | F[T_n \mid \mathcal{F}_S] - L_S | \leq  \frac{1}{n} \enspace .
 \end{equation}
 Once again, we can assume $U_n \leq U_1$. Define
 \begin{equation}
 \label{equation:def-Vnagain}
 V_n ~  := ~ \left\{\omega\in \Omega \, \left| \, \frac{1}{n} \leq \RVeps(\omega) <  \frac{1}{n-1} \right. \right\} \enspace ,
 \end{equation}
 where in the case $n = 1$ we understand $\frac{1}{n-1}$ as $\infty$. Clearly $\{V_n\}_n$ is an $\mathcal{F}_S$-partition of $\Omega$. 
 Define a stopping time by $U := \sum_n \, \One_{V_n} \cdot U_n$.  Clearly $S < U \leq U_1$.
 We show that $U$ satisfies the property required to establish   randomly uniform convergence \ref{s-l:a}. Suppose that $T$ is such that $S < T \leq U$. We need to show that $ \Prob[\, | F[T \mid \mathcal{F}_S] - L_S | \leq  \RVeps] = 1$.

Define
 \[
 T_n(\omega) ~ := ~ \begin{cases}
   T(\omega) & \text{if $\omega \in V_n$} \\
    U_n(\omega) & \text{if $\omega \not\in V_n$}  \enspace .
   \end{cases}
 \]
 Because $T \leq U$, we have $T_n \leq U_n$. Whence, since $\{V_n\}_n$ is an $\mathcal{F}_S$-partition of $\Omega$, we have
  \begin{align*}
  & \Prob[ | F[T \mid \mathcal{F}_S] - L_S | \leq  \RVeps]  \\
  & ~~~ = ~  \Prob\left[ \, \left|\,  F\left[ \left. \sum_n \One_{V_n} \cdot T_n \, \right| \,  \mathcal{F}_S\right] - L_S \right| \leq  \RVeps \right]  \\[1ex]
  & ~~~ = ~ \Prob\left[ \, \left|\, \left(\sum_n \One_{V_n} \cdot F[T_n \mid \mathcal{F}_S] \right)  - L_S \right| \leq  \RVeps \right] & & \text{by the partition property}
  \\[1ex]
  & ~~~ = ~ \sum_n\, \Prob(\{\omega \in V_n \mid \, |F[T_n \mid \mathcal{F}_S](\omega) - L_S(\omega)| \leq \RVeps(\omega) \})
  \\
  & ~~~ = ~ \sum_n\, \Prob(V_n) & & \text{by~\eqref{equation:one-on-n} and \eqref{equation:def-Vnagain}}
  \\
  & ~~~ = ~ 1 & & \text{because $\Omega  = \biguplus_{i = 1}^\infty V_i$ \enspace .}
  \end{align*}
Again both sides in the partition property are defined because all stopping times are $\leq U_1$.
 \end{proof}
 
%
%

 \begin{definition}[Stopping limit] 
 \label{definition:stopping-limit}
 An $\mathcal{F}_S$-measurable random variable $L_S$ is a \emph{stopping limit} for a stopping functional $F$ at a finite stopping time $S$ (notation $L_S = \lim_{T \downarrow S} \,F[T \mid \mathcal{F}_S]$)
if any of the equivalent properties of Definition~\ref{definition:net-convergence}  hold.

A right-continuous process $(L)$ is a \emph{stopping-limit process} for  $F$ if, for every 
finite stopping time $S$, it holds that $L_S$ is a stopping limit for $F$ at $S$.  \end{definition}
 
 \noindent
 The notion of stopping-limit process could be defined more generally for a progressive process $(L)$, but the assumption of  right continuity seems general enough for practical purposes and ties in with a general ethos underlying the present work, in which time is viewed in a right-oriented way. 
 
It is immediate from the definition that if $L_S$ is a stopping limit for $F$ at $S$ and   $L'_S$ is $\mathcal{F}_S$-almost surely equal to $L_S$ then $L'_S$ is also a stopping limit for  $F$ at $S$. 
Similarly, if $(L)$ is a stopping-limit process for $F$ and $(L')$ is indistinguishable from $(L)$ then $(L')$ is also a stopping-limit process for $F$. 
While the result below is basic, the proof does illustrate one subtle point:  in the four
equivalent properties of Definition~\ref{definition:net-convergence}, it is important that we quantify over stopping times $T$ with $T \leq U$ rather than over stopping times $T$ with $T < U$. Indeed, the stopping time $T := U \wedge U'$ used in the proof below need not satisfy $T < U$. Moreover, there seems to be no general possibility of resolving this by by finding an intermediate stopping time  $T$ such that $S < T < U \wedge U'$, since  such an interpolating stopping time need not exist. (For example, if $U$ is the first jump time of a Poisson process then $0 < U$, but there does not exist any stopping time $T$ with $0 < T < U$.)

\begin{proposition}
\label{proposition:limit-unique}
\begin{enumerate}
\item If a stopping functional $F$ has a stopping limit  at a stopping time  $S< \infty$  then this is uniquely determined up to $\mathcal{F_S}$-almost-sure equality. 

\item If  a stopping-limit process for $F$ exists then  it is uniquely determined up to indistinguishability.

\end{enumerate}
\end{proposition}

\begin{proof}
For statement 1, suppose that $L_S$ and $L'_S$ are two $\mathcal{F}_S$-measurable random variables that are not $\mathcal{F}_S$-almost-surely equal.
Then  there exist a positive probability event $W \in \mathcal{F}_S$ and $\delta > 0$ such that, for all $\omega \in W$, it holds that $|L'_S(\omega) - L_S(\omega)| > \delta$. 
Assume, for contradiction, that 
both $L_S$ and $L'_S$ are stopping limits for $F$ at $S$. 
Using $\varepsilon := \frac{\delta}{2}$ and the definition of uniform limits, let $U$ and $U'$ be finite stopping times such that: for all $T$ with $S \leq T \leq U$, we have
\[ |F[T \mid \mathcal{F}_S] - L_S | \leq \varepsilon \enspace ,\]
and, for all $T$ with $S \leq T \leq U'$, we have
\[ |F[T \mid \mathcal{F}_S] - L'_S | \leq \varepsilon \enspace .\]
Since both inequalities are true for  $T := U \wedge U'$, we have  $|L'_S - L_S| \leq 2\varepsilon$. But this contradicts the fact that 
$|L'_S(\omega) - L_S(\omega)| > \delta$ for all $\omega \in W$. 

Statement 2 follows since, by 1, any two stopping-limit processes are modifications of one another, and  it is standard that any two right-continuous processes that are modifications of each other are indistinguishable.
\end{proof}


The notion of stopping limit induces a natural notion of continuity for stopping functionals.
\begin{definition}[Stopping continuity]
\label{definition:stopping-continuous}
We say that a stopping functional $F$ is \emph{stopping continuous at} a finite stopping time $S$ if 
\begin{equation}
\label{equation:stopping-continuous}
F[S\mid \mathcal{F}_S] = \lim_{T \downarrow S} \, F[T \mid \mathcal{F}_S] \enspace .
\end{equation}
The stopping functional $F$ itself is said to be \emph{stopping continuous} if it is stopping continuous at every finite stopping time $S$.
\end{definition}

One of the consequences of our use of generalised conditional expectation is that, 
if $(X)$ is a continuous process, then the stopping functionals for conditional expectation
$\Exp[T \mid \mathcal{F}_S]$ and conditional variance $\Var[T \mid \mathcal{F}_S]$ are always defined
for sufficiently small stopping times $T > S$. 
In fact, both stopping functionals are even stopping continuous.
This follows from Propositions~\ref{proposition:exp-continuous} and~\ref{proposition:var-continuous} below, which establish that stopping continuity holds more generally for cadlag processes, as long as we replace $X_T$ in~\eqref{def:exp} and~\eqref{def:var} with $X_{T-}$ as defined in~\eqref{eqn:T-minus} in Section~\ref{section:preliminaries}. That is, we prove continuity for the \emph{cadlag-adjusted} stopping functionals
for conditional expectation and conditional variance, $\Exp[X_{T-} \mid \mathcal{F}_S]$  and $\Var[X_{T-} \mid \mathcal{F}_S]$ respectively. Note that, in the special case of a continuous process, the cadlag-adjusted stopping functionals coincide with the original ones. 

\begin{proposition}
\label{proposition:exp-continuous}
If $(X)$ is a cadlag process then the cadlag-adjusted stopping functional for conditional expectation
$\Exp[X_{T-} \mid \mathcal{F}_S]$  is stopping continuous. 
\end{proposition}
\begin{proof} We show
\begin{equation}
\label{equation:exp-continuous-equalities}
\lim_{T \downarrow S} \, \Exp[X_{T-} \mid \mathcal{F}_S] ~  =  ~ X_S  ~ =  ~ \Exp[X_{S-} \mid \mathcal{F}_S]  \enspace .
\end{equation}
The right-hand equality holds because~\eqref{eqn:T-minus} defines $X_{S-} := X_S$.
We prove the left-hand equality using uniform convergence. Accordingly, let $\varepsilon > 0$ be arbitrary. Define a first-approach stopping time, using the notation  introduced in~\eqref{equation:plus-minus},
\[
U(\omega) ~ := ~ (S(\omega) + 1) \wedge \inf \{t >  S(\omega) \mid \, |X_{t\pm}(\omega) - X_S(\omega) | \geq \epsilon \} \enspace ,
\]
where the role of $S + 1$ is to ensure the  finiteness of U. Clearly $U \geq S$. To prove $U > S$, consider any $\omega \in \Omega$. By the right-continuity of $t \mapsto X_t(\omega)$, there exists $\delta > 0$ such that, for all
$t$ with $S(\omega) \leq t < S(\omega) + \delta$,  it holds that  $| X_t(\omega) - X_S(\omega) | < \frac{\varepsilon}{2}$.
It follows that $\inf \{t >  S(\omega) \mid \, |X_{t\pm}(\omega) - X_S(\omega) | \geq \epsilon \} \geq S(\omega) + \delta$,
whence $U(\omega) > S(\omega)$.

To establish the left-hand equality of \eqref{equation:exp-continuous-equalities}, 
let $T$ be such that $S \leq T \leq U$. 
Because $T \leq U$, we have, by the definition of $U$, that $|X_{T-}(\omega)  - X_S(\omega)| \leq \varepsilon$ for every $\omega \in \Omega$ (equality can hold in the case $T(\omega) = U(\omega)$). Hence $|\Exp[X_{T-} \mid \mathcal{F}_S] - X_S| \leq \varepsilon$, which is what we were required to show. 
\end{proof}

\begin{proposition}
\label{proposition:var-continuous}
If $(X)$ is a cadlag process then the cadlag-adjusted stopping functional for conditional variance
\begin{equation}
\label{equation:variance-rate}
\Var[X_{T-} | \mathcal{F}_S] ~ := ~ \Exp[(X_{T-}  - \Exp[X_{T-} | \mathcal{F}_S])^2 \,| \mathcal{F}_S]\, \enspace .
\end{equation}
 is stopping continuous. 
\end{proposition}
\begin{proof} We show
\begin{equation}
\label{equation:var-continuous-equalities}
\lim_{T \downarrow S} \, \Var[X_{T-} \mid \mathcal{F}_S] ~  =  ~ 0  ~ =  ~ \Var[X_{S-} \mid \mathcal{F}_S]  \enspace .
\end{equation}
Again, the right-hand equality holds because $X_{S-} = X_S$, and 
we prove the left-hand equality using uniform convergence. Let $\varepsilon > 0$ be arbitrary. 

Using the stopping continuity of cadlag-adjusted conditional expectation, let $U_1 > S$ be a finite stopping time such that,
for all $T$ with $S \leq T \leq U_1$, we have $|\Exp[X_{T-} \mid \mathcal{F}_S] - X_S| \leq \frac{\sqrt{\varepsilon}}{2}$.
Define 
\[
U_2(\omega) ~ := ~ (S(\omega) + 1) \wedge \inf \left\{t >  S(\omega)\, \left|  \, |X_{t\pm}(\omega) - X_S(\omega) | \geq  \frac{\sqrt{\varepsilon}}{2} \right.\right\} \enspace .
\]
Then $U_2 > S$, as in the proof of Proposition~\ref{proposition:exp-continuous}.
Define $U := U_1 \wedge U_2$. Clearly $U > S$.

To establish the left-hand equality of \eqref{equation:var-continuous-equalities}, 
let $T$ be such that $S \leq T \leq U$. 
Because $T \leq U_1$, we have $|\Exp[X_{T-} \mid \mathcal{F}_S] - X_S| \leq \frac{\sqrt{\varepsilon}}{2}$.
Because $T \leq U_2$, we have $|X_{T-}  - X_S| \leq \frac{\sqrt{\varepsilon}}{2}$. 
Hence $(X_{T-} - \Exp[X_{T-} \mid \mathcal{F}_S])^2 \leq \varepsilon$. It follows that
$\Exp[(X_{T-} - \Exp[X_{T-} \mid \mathcal{F}_S])^2  \mid \mathcal{F}_S] \leq \varepsilon$, 
which is what we needed to show. 
\end{proof}



\section{Stopping derivatives}
\label{section:stopping-derivatives}

 Using the notion of stopping limit from Definition~\ref{definition:stopping-limit}, the definition of  stopping derivative takes the expected form for a right derivative,  cf.~\eqref{equation:drift-as-sd}.
  \begin{definition}[Stopping derivative]
 \label{definition:stopping-derivative}
An $\mathcal{F}_S$-measurable $\Real$-valued random variable $D_S$ is a
 \emph{stopping derivative} for  a stopping functional $F$ at a stopping time $S < \infty$  if
 \begin{equation}
 \label{equation:stopping-derivative}
 D_S ~ = ~ \lim_{T \downarrow S} \, \frac{F[T \mid \mathcal{F}_S] - F[S \mid \mathcal{F}_S]}{\Exp[T-S \mid \mathcal{F}_S]} \enspace . 
 \end{equation}
 \noindent
 A right-continuous process $(D)$ is a \emph{stopping-derivative process} for  $F$ if, for every 
finite stopping time $S$, it holds that $D_S$ is a stopping derivative for $F$ at $S$.
\end{definition}

Because stopping derivatives are defined as stopping limits, they inherit the uniqueness properties of stopping limits, as given by 
Proposition~\ref{proposition:limit-unique}. That is, the stopping derivative at $S$ is unique up to $\mathcal{F_S}$-almost-sure equality. Also, the stopping-derivative process for $F$, if it exists, is determined up to indistinguishability.
Adapting standard derivative notation, we write
$F'[\mathcal{F}_S]$ for the stopping derivative of $F$ at $S$, if it exists, and
$(F'[\mathcal{F}_t]) _{t \geq 0}$ (more concisely $(F'[\mathcal{F}])\,$) for 
a stopping-derivative process for $F$.

In real analysis,  continuity is a necessary condition for differentiability. Similarly, stopping continuity is a necessary condition for stopping derivatives to exist.
\begin{proposition}
\label{proposition:derivative-continuous}
If a stopping functional $F$ has a stopping derivative at a finite stopping time $S$ then $F$ is stopping continuous at $S$.
\end{proposition}
 \begin{proof} Let $D_S$ be a stopping derivative for  $F$ at $S$. We prove~\eqref{equation:stopping-continuous}, using
 uniform convergence. Suppose $\varepsilon > 0$. 
 
 Let $U_1 >S$ be a finite stopping time, as guaranteed by $D_S$ being a stopping derivative, such that, for all $T$ with $S < T \leq U_1$, it holds that
 \begin{equation}
 \label{equation:dc-frac}
 \left| \frac{F[T \mid \mathcal{F}_S] - F[S \mid \mathcal{F}_S]}{\Exp[T-S \mid \mathcal{F}_S]} - D_S \right| ~ \leq ~ 1 \enspace .
 \end{equation}
 Define 
 \[ U ~  := ~ U_1 \wedge \left(S + \frac{\varepsilon}{|D_S|+1}\right)\enspace .\]
 Clearly $U > S$. 
 
 To establish the stopping continuity of $F$ at $S$, let $T$ be such that $S < T \leq U$. We need to prove 
 \begin{equation}
 \label{equation:dc-toshow}
 |F[T \mid \mathcal{F}_S] - F[S \mid \mathcal{F}_S]| \leq \varepsilon \enspace .
  \end{equation}
 Because $T \leq U_1$, it follows from~\eqref{equation:dc-frac} that
 \begin{equation}
 \label{equation:derivative-continuous}
 |F[T \mid \mathcal{F}_S] - F[S \mid \mathcal{F}_S]| ~ \leq ~ (|D_S| + 1) \Exp[T-S \mid \mathcal{F}_S] \enspace .
 \end{equation}
 Also because $T \leq U$, we  have 
 \begin{equation}
 \label{equation:dc-last}
 \Exp[T-S \mid \mathcal{F}_S] ~ \leq ~ \frac{\varepsilon}{|D_S|+1} \enspace . 
 \end{equation}
 Clearly~\eqref{equation:dc-toshow} follows from~\eqref{equation:derivative-continuous} and~\eqref{equation:dc-last}.
  \end{proof}

In the sequel, we shall have recourse to consider many instantiations of the notion of stopping derivative. 
At this point, we define just the two principal instances of interest, the \emph{drift} and \emph{variance rate}, as
motivated in Section~\ref{section:intro}.

\begin{definition}[Drift]
\label{definition:drift}
The \emph{drift} of a cadlag process $(X)$ at $S$ is a stopping derivative for the \emph{cadlag-adjusted conditional expectation} function $\Exp[X_{T-} | \mathcal{F}_S]\,$.
\end{definition}

\begin{definition}[Variance rate]
\label{definition:variance-rate}
The \emph{variance rate} of a cadlag process $(X)$ at $S$ is a stopping derivative for the \emph{cadlag-adjusted conditional variance} function $\Var[X_{T-} | \mathcal{F}_S]$, defined in~\eqref{equation:variance-rate}.
\end{definition}
\noindent
In the above definitions, the cadlag-adjusted versions of conditional expectation and variance are needed to obtain definitions that are applicable to cadlag processes, since
it is the cadlag-adjusted versions that are stopping continuous (Propositions~\ref{proposition:exp-continuous} and~\ref{proposition:var-continuous}).


The concepts of drift and variance rate will be studied in detail in Sections~\ref{section:drift} and~\ref{section:diffusion} respectively. To set up the necessary technical machinery, we develop, in the remainder of the present section, some of the general theory of stopping derivatives. 

The first result establishes a useful characterisation of stopping derivatives. 
\begin{proposition}[Characterisation of stopping derivatives]
\label{proposition:csd}
An $\mathcal{F}_S$-measurable $\Real$-valued random variable $D_S$ is a
 stopping derivative for  a stopping functional $F$ at a stopping time $S < \infty$  if and only if, for every 
 $\mathcal{F}_S$-measurable random variable $\underline{\varepsilon} > 0$, there exists a stopping time
 $U> S$ such that, for all finite stopping times $T$ with $S \leq T \leq U$,
 \begin{equation}
 \label{equation:csd}
(D_S - \underline{\varepsilon}\,) \,  \Exp[T-S  \mid \mathcal{F}_S] ~ \leq ~ 
F[T \mid \mathcal{F}_S] - F[S \mid \mathcal{F}_S] ~ \leq ~ (D_S + \underline{\varepsilon}\,) \, \Exp[T-S  \mid \mathcal{F}_S] \enspace .
 \end{equation}
 \end{proposition}

\noindent
One advantage of the fraction-free formulation of stopping derivative in  the proposition is that it asks 
only for  the non-strict inequality $T \geq S$ instead of a strict inequality $T > S$. This is helpful when working with stopping times, because, unlike with deterministic time, there is no simple dichotomy between just two cases $T > S$ (almost surely) and $T = S$ (almost surely), and it  sometimes proves technically convenient
to work with stopping times $T$, for which both $T=S$ and $T>S$ hold on positive measure subsets.  (This is helpful, for example,  in the proof of Theorem~\ref{theorem:driftless-martingale}.) The other (relatively trivial) change is that there is no  requirement in the proposition  for  the stopping time $U$ be finite. 

\begin{proof}[Proof of Proposition~\ref{proposition:csd}]
We prove that the characterisation in Proposition~\ref{proposition:csd} is equivalent to Definition~\ref{definition:stopping-derivative} with the latter being  understood in terms of randomly uniform convergence. 

Firstly, assume that the characterisation in Proposition~\ref{proposition:csd} holds. To show~\eqref{equation:stopping-derivative},
let $\RVeps > 0$ be $\mathcal{F}_S$-measurable. Let $U > S$ be such that, for any finite $T$ with $S \leq T \leq U$,~\eqref{equation:csd} holds. Define
$U' = U \wedge (S + 1)$. Then $U'$ is a finite stopping time with $S < U' \leq U$. Let $T$ be such that $S < T \leq U'$.
Because $\Exp[T - S \mid \mathcal{F}_S] > 0$, It is immediate from~\eqref{equation:csd} that
\begin{equation}
\label{equation:proof:csd}
\left|\, \frac{F[T \mid \mathcal{F}_S] - F[S \mid \mathcal{F}_S]}{\Exp[T-S \mid \mathcal{F}_S]}  - D_S\, \right| \leq \RVeps \enspace .
\end{equation}

Conversely, assume that~\eqref{equation:stopping-derivative} holds. Let $\RVeps > 0$ be $\mathcal{F}_S$-measurable.
Applying~\eqref{equation:stopping-derivative}, let $U > S$ be a finite stopping time such that, for any $T$ with $S < T \leq U$,~\eqref{equation:proof:csd} holds. Let $T$ be such that  $S \leq T \leq U$.
We need to show that~\eqref{equation:csd} holds. Define 
\[
V ~ := ~ \{ \omega \in \Omega \mid T(\omega) > S(\omega)\} \enspace .
\]
Clearly $V \in \mathcal{F}_S$. Define a stopping time $T_V := \One_V \cdot T + \One_{V^C} \cdot U$. Clearly $S < T_V \leq U$.
Hence,~\eqref{equation:proof:csd} holds for $T_V$, i.e.,
\[
\left|\, \frac{F[\One_V \cdot T + \One_{V^C} \cdot U \mid \mathcal{F}_S] - F[\One_V \cdot S  + \One_{V^C} \cdot S \mid \mathcal{F}_S]}{\Exp[(\One_V \cdot T + \One_{V^C} \cdot U)-(\One_V \cdot S  + \One_{V^C} \cdot S) \mid \mathcal{F}_S]}  - D_S\, \right| \leq \RVeps \enspace .
\]
By the partition property, which applies because all stopping times are $\leq U$, and by laws of conditional expectation,
\[
\left|\, \frac{\One_V \cdot (F[T  \mid \mathcal{F}_S] - F[S  \mid \mathcal{F}_S]) + 
\One_{V^C} \cdot (F[U  \mid \mathcal{F}_S] - F[S  \mid \mathcal{F}_S])}
{\One_V \cdot \Exp[T -S \mid \mathcal{F}_S] + \One_{V^C} \cdot  \Exp[U -S \mid \mathcal{F}_S]} - D_S\, \right| \leq \RVeps \enspace .
\]
So for almost all $\omega \in V$, we have
\[
\left|\, \frac{F[T  \mid \mathcal{F}_S](\omega)- F[S  \mid \mathcal{F}_S](\omega)}
{\Exp[T -S \mid \mathcal{F}_S](\omega)} - D_S(\omega)\, \right| \leq \RVeps(\omega) \enspace ,
\]
whence
\[
((D_S - \underline{\varepsilon})  \Exp[T-S  \mid \mathcal{F}_S])(\omega)~ \leq ~ 
(F[T \mid \mathcal{F}_S] - F[S \mid \mathcal{F}_S])(\omega) ~ \leq ~ ((D_S+ \underline{\varepsilon})  \Exp[T-S  \mid \mathcal{F}_S]) (\omega) .
\]
Moreover, for almost all $\omega \in V^C$, we have $T(\omega)= S(\omega)$, so
\[
((D_S - \underline{\varepsilon})  \Exp[T-S  \mid \mathcal{F}_S])(\omega)= 
(F[T \mid \mathcal{F}_S] - F[S \mid \mathcal{F}_S])(\omega) =  ((D_S+ \underline{\varepsilon})  \Exp[T-S  \mid \mathcal{F}_S]) (\omega) = 0.
\]
The two cases above, for  $\omega \in V$ and  $\omega \in V^C$, combine to give us~\eqref{equation:csd}, as required.
\end{proof}

We devote the remainder of the section to developing the basic calculus of stopping derivatives. Throughout the development, we use~\eqref{equation:csd} above as our working definition of stopping derivative. We remark that, by working directly with this characterisation of stopping derivatives, we never again in this section use  the partition property of stopping functionals.

Stopping derivatives enjoy analogues of the usual rules of differential calculus. We present a selection of such properties: linearity, product and chain rules and a time-change rule. 
The proofs are fundamentally straightforward---indeed, they are very similar to the usual ones. Nevertheless, the use of stopping times together with a random $\underline{\varepsilon}$ adds some technical overhead to the arguments, making the the proofs that follow somewhat cumbersome. Readers who are keen to arrive at  the development of stochastic calculus based on drift and variance rate may prefer to  skip the remainder of present section and refer back to it only when needed.

Given a stopping functional $F$ and a random variable $A$, we write $AF$ for the stopping functional
\[
S \leq T ~ \mapsto ~ \begin{cases}
AF[T \mid \mathcal{F}_S] & \text{if $A$ is $\mathcal{F}_S$-measurable}
\\
\text{undefined} & \text{otherwise.}
\end{cases}
\]

\begin{proposition}[Linearity]
\label{proposition:linear}
If $F$ and $G$ have stopping derivatives at a finite stopping time $S$ and $A,B$ are 
$\mathcal{F}_S$-measurable random variables, then $A F + B G$ has a stopping derivative at $S$ and
\[
(A F + B G)'\,[\mathcal{F}_S] ~ = ~ A F'[\mathcal{F}_S]  + B G'[\mathcal{F}_S] \enspace .
\]
\end{proposition}
\noindent
We omit the straightforward proof.
%
%

\begin{proposition}[Product rule]
\label{proposition:product}
If $F$ and $G$ have stopping derivatives at a finite stopping time $S$ then so does $F G$ and:
\[
(FG)'\,[\mathcal{F}_S] ~ = ~ F[S \mid \mathcal{F}_S]\, G'[\mathcal{F}_S]  + G[S \mid \mathcal{F}_S]\,F'[\mathcal{F}_S] \enspace .
\]
\end{proposition}

\begin{proof}
Consider any $\mathcal{F}_S$-measurable $\underline{\varepsilon} > 0$. Define
\[
\underline{\varepsilon}_1 ~ := ~ \frac{\underline{\varepsilon}}{3\,(|G[S \mid \mathcal{F}_S]| \vee 1)}
\qquad
\underline{\varepsilon}_2 ~ := ~ \frac{\underline{\varepsilon}}{3\,(|F[S \mid \mathcal{F}_S]| \vee 1)} \enspace .
\]
Let $U_1 > S$ be such that, 
for any finite $T$ with $S \leq  T \leq U_1$,
\begin{equation}
\label{equation:product-a}
\left(F'[\mathcal{F}_S]  - \underline{\varepsilon}_1\right) \Exp[T-S \mid \mathcal{F}_S] ~ \leq ~ 
F[T \mid \mathcal{F}_S] - F[S \mid \mathcal{F}_S]
~ \leq ~\left(F'[\mathcal{F}_S]  + \underline{\varepsilon}_1\right) \Exp[T-S \mid \mathcal{F}_S]
 \enspace ,
\end{equation}
and $U_2 > S$ be such that, 
for any finite $T$ with $S \leq  T \leq U_2$,
\begin{equation}
\label{equation:product-b}
\left(G'[\mathcal{F}_S]  - \underline{\varepsilon}_2 \right) \Exp[T-S \mid \mathcal{F}_S] ~ \leq ~ 
G[T \mid \mathcal{F}_S] - G[S \mid \mathcal{F}_S]
~ \leq ~\left(G'[\mathcal{F}_S]  + \underline{\varepsilon}_2\right)\Exp[T-S \mid \mathcal{F}_S]
 \enspace ,
\end{equation}
Define 
\[U_3 ~ := ~ S + \frac{\underline{\varepsilon}}{3\, ((|F'[\mathcal{F}_S]| + \underline{\varepsilon}_1) (|G'[\mathcal{F}_S]| + \underline{\varepsilon}_2) \vee 1)} \enspace .
\]
Clearly $U_3$ is a stopping time with $S  < U_3$. Define $U := U_1 \wedge U_2 \wedge U_3$, and consider any 
finite $T$
with $S \leq  T \leq U$. 

Multiplying~\eqref{equation:product-a} by $G[S \mid \mathcal{F}_S]$, we obtain
\begin{align}
\nonumber
&  \left(G[S \mid \mathcal{F}_S]\, F'[\mathcal{F}_S]  - D_1 \right) \Exp[T-S \mid \mathcal{F}_S] 
\\ 
\nonumber
& \qquad \leq ~ 
F[T \mid \mathcal{F}_S]\, G[S \mid \mathcal{F}_S] - F[S \mid \mathcal{F}_S]\, G[S \mid \mathcal{F}_S] \\
\label{equation:product-f}
& \qquad \qquad  \leq ~\left(G[S \mid \mathcal{F}_S]\, F'[\mathcal{F}_S]  + D_1\right) \Exp[T-S \mid \mathcal{F}_S]
 \enspace ,
\end{align}
where $D_1 := |G[S \mid \mathcal{F}_S]|\, \underline{\varepsilon}_1$. 
Similarly, multiplying~\eqref{equation:product-b} by $F[T \mid \mathcal{F}_S]$, we obtain
\begin{align}
\nonumber
&  \left(F[T \mid \mathcal{F}_S]\, G'[\mathcal{F}_S]  - |F[T \mid \mathcal{F}_S]|\, \underline{\varepsilon}_2 \right) \Exp[T-S \mid \mathcal{F}_S] 
\\ 
\nonumber
& \qquad \leq ~ 
F[T \mid \mathcal{F}_S]\, G[T \mid \mathcal{F}_S] - F[T \mid \mathcal{F}_S]\, G[S \mid \mathcal{F}_S] 
\\
\label{equation:product-c}
& \qquad \qquad  \leq ~\left(F[T \mid \mathcal{F}_S]\, G'[\mathcal{F}_S]  + |F[T \mid \mathcal{F}_S]|\, \underline{\varepsilon}_2 \right) \Exp[T-S \mid \mathcal{F}_S]
 \enspace .
\end{align}
Then  using~\eqref{equation:product-a} to replace $F[T \mid \mathcal{F}_S]$ with $F[S \mid \mathcal{F}_S]$ in the bounds of~\eqref{equation:product-c},
  we obtain
\begin{align}
\nonumber
&  \left(F[S \mid \mathcal{F}_S]\, G'[\mathcal{F}_S]  - D_2 \right) \Exp[T-S \mid \mathcal{F}_S] 
\\ 
\nonumber
& \qquad \leq ~ 
F[T \mid \mathcal{F}_S]\, G[T \mid \mathcal{F}_S] - F[T \mid \mathcal{F}_S]\, G[S \mid \mathcal{F}_S] 
\\
\label{equation:product-e}
& \qquad \qquad  \leq ~\left(F[S \mid \mathcal{F}_S]\, G'[\mathcal{F}_S]  + D_2  \right) \Exp[T-S \mid \mathcal{F}_S]
 \enspace ,
\end{align}
where 
\[D_2  ~ :=  ~(|F'[\mathcal{F}_S]| + \underline{\varepsilon}_1) \, (|G'[\mathcal{F}_S]| + \underline{\varepsilon}_2) \,\Exp[T-S \mid \mathcal{F}_S] 
+ |F[S \mid \mathcal{F}_S]|\, \underline{\varepsilon}_2 \enspace .\]
Adding~\eqref{equation:product-e} and~\eqref{equation:product-f}, we get:
\begin{align*}
\nonumber
&  \left(F[S \mid \mathcal{F}_S]\, G'[\mathcal{F}_S]  +  G[S \mid  \mathcal{F}_S]\, F'[\mathcal{F}_S] - D_1 - D_2
 \right) \Exp[T - S \mid \mathcal{F}_S] 
\\ 
\nonumber
& \qquad \leq ~ 
F[T \mid \mathcal{F}_S]\, G[T \mid \mathcal{F}_S] - F[S \mid \mathcal{F}_S]\, G[S \mid \mathcal{F}_S] 
\\
& \qquad \qquad  \leq ~
\left(F[S \mid \mathcal{F}_S]\, G'[\mathcal{F}_S]  +  G[S \mid  \mathcal{F}_S]\, F'[\mathcal{F}_S] + D_1 + D_2
 \right) \Exp[T - S \mid \mathcal{F}_S] \enspace .
\end{align*}
But $D_1 \leq \frac{\underline{\varepsilon}}{3}$, by the definition of $\underline{\varepsilon}_1$, and
$D_2 \leq \frac{2\, \underline{\varepsilon}}{3}$, by the definitions of $\underline{\varepsilon}_2$ and $U_3$. Hence
\begin{align*}
\nonumber
&  \left(F[S \mid \mathcal{F}_S]\, G'[\mathcal{F}_S]  +  G[S \mid  \mathcal{F}_S]\, F'[\mathcal{F}_S] - \underline{\varepsilon}
 \right) \Exp[T - S \mid \mathcal{F}_S] 
\\ 
\nonumber
& \qquad \leq ~ 
F[T \mid \mathcal{F}_S]\, G[T \mid \mathcal{F}_S] - F[S \mid \mathcal{F}_S]\, G[S \mid \mathcal{F}_S] 
\\
& \qquad \qquad  \leq ~
\left(F[S \mid \mathcal{F}_S]\, G'[\mathcal{F}_S]  +  G[S \mid  \mathcal{F}_S]\, F'[\mathcal{F}_S] + \underline{\varepsilon}
 \right) \Exp[T - S \mid \mathcal{F}_S] \enspace .
\end{align*}
\end{proof}

To obtain a chain rule with an appropriate level of generality, we consider the composition of $F$ with a \emph{random} function.
\begin{definition}[Jointly measurable random function]
\label{definition:random-function}
Given a sub-$\sigma$-algebra $\mathcal{G} \subseteq \mathcal{E}$, we say that a function 
$\underline{f} \colon \Omega \times \Real \to \Real$
is a \emph{$\mathcal{G}$-jointly-measurable  random function} from $\Real$ to $\Real$ if it is 
measurable from the product $\sigma$-algebra $\mathcal{G} \otimes \mathcal{B}$ to $\mathcal{B}$.
(In general, we allow
$\underline{f}$ to be partial, with a $(\mathcal{G} \otimes \mathcal{B})$-measurable subset of $\Omega \times \Real$ as its domain of definedness.)  
\end{definition}

When working with random functions, we write $\underline{f}_\omega$ for the function $x \mapsto \underline{f}(\omega,x) : \Real \to \Real$. 
We also require  measurable versions of the notions of continuity and differentiability.
\begin{definition}[Measurable continuity and differentiability]
\label{definition:measurable-continuity-differentiability}
 Suppose $\underline{f}$  is 
 $\mathcal{G}$-jointly-measurable  and $X$ is a $\mathcal{G}$-measurable random variable.
 \begin{enumerate}
 \item 
 $\underline{f}$ is \emph{$\mathcal{G}$-measurably continuous} at  $X$ if, for every $\mathcal{G}$-measurable
$\RVeps > 0$, there exists a $\mathcal{G}$-measurable
$\RVdel > 0$ such that, for  $\mathcal{G}$-almost-all  $\omega$,
\[
\forall y, ~~  |y - X(\omega)| \leq \RVdel(\omega) \text{~implies~}
\left| \underline{f}_\omega(y) - \underline{f}_\omega(X(\omega)) \right| ~ \leq ~ \RVeps(\omega) \enspace ,
\]
where we are implicitly assuming that $\underline{f}_\omega(y)$ and $\underline{f}_\omega(X(\omega))$ are defined.
\item $\underline{f}$ is \emph{$\mathcal{G}$-measurably differentiable} at $X$ if, 
there exists a $\mathcal{G}$-measurable random variable $\underline{f}'(X)$ such that, for every $\mathcal{G}$-measurable
$\RVeps > 0$, there exists a $\mathcal{G}$-measurable
$\RVdel > 0$ such that, for  $\mathcal{G}$-almost-all  $\omega$,
\[
\forall y, ~~ 0 <  |y - X(\omega)| \leq \RVdel(\omega) \text{~implies~} 
\left| \frac{\underline{f}_\omega(y) - \underline{f}_\omega(X(\omega))}{y - X(\omega)} - 
\underline{f}'(X) (\omega) \right| ~ \leq ~ \RVeps(\omega) \enspace ,
\]
where we are again assuming that $\underline{f}_\omega(y)$ and $\underline{f}_\omega(X(\omega))$ are defined.
\end{enumerate}
 \end{definition}
 \noindent
 (The use of non-strict inequalities in 1 and 2 above is for mild technical convenience in proofs. Of course the definitions remain equivalent if all inequalities are strictified.)
 
 A useful sufficient condition for  
$\mathcal{G}$-measurable continuity at $X$ to hold is the following: there exists a measure $1$ set $W \in \mathcal{G}$ and a
$\mathcal{G} \otimes \mathcal{B}$-measurable set $D \subseteq W \times \Real$, such that, for every $\omega \in W$,
the set $D_\omega := \{x \in \Real \mid (\omega,x ) \in D\}$ is an open interval containing $X(\omega)$ and the function
 $\underline{f}_\omega$ is defined and continuous on $D_\omega$. In the case that $\underline{f}_\omega$ is also 
 differentiable at $X(\omega)$ for every $\omega \in W$, it further holds that  $\underline{f}$ is $\mathcal{G}$-measurably differentiable.
 These observations motivate the following definition of what it means for a random function to be measurably continuously differentiable,
 which will be used later in our formulation of the It\^o formula (Theorem~\ref{theorem:ito-formula}).

\begin{definition}[Measurable continuous differentiability]
\label{definition:measurable-continuous-differentiability}
Suppose $\underline{f}$  is 
 $\mathcal{G}$-jointly-measurable  and $X$ is a $\mathcal{G}$-measurable random variable.
We say that $\underline{f}$ is \emph{$\mathcal{G}$-measurably continuously differentiable} at  $X$ if
there exists a measure $1$ set $W \in \mathcal{G}$ and a
$\mathcal{G} \otimes \mathcal{B}$-measurable set $D \subseteq W \times \Real$, such that, for every $\omega \in W$,
the set $D_\omega := \{x \in \Real \mid (\omega,x ) \in D\}$ is an open interval containing $X(\omega)$ and the function
 $\underline{f}_\omega$ is defined and continuously differentiable on $D_\omega$.
\end{definition}

Is is routine to show that if $\underline{f}$ is \emph{$\mathcal{G}$-measurably continuously differentiable} at  $X$
with respect to $D \subseteq W \times \Real$ as above, then the derivative 
\[\omega, x \mapsto \underline{f}_\omega'(x) : D \to \Real\]
is itself a $\mathcal{G}$-measurably continuous random function.

Given a stopping functional $F$ and a jointly-measurable random function $\underline{f}$, the composite stopping functional
$\underline{f} \circ F$ is defined by
\[
S \leq T ~ \mapsto ~ \begin{cases}
\underline{f} (F[T \mid \mathcal{F}_S]) & \text{if $\underline{f}$ is $\mathcal{F}_S$-jointly-measurable}
\\
\text{undefined} & \text{otherwise.}
\end{cases}
\]

\begin{proposition}[Chain rule]
Suppose $F$ has a stopping derivative at a finite stopping time $S$, and
$\underline{f}$ is an $\mathcal{F}_S$-jointly-measurable random function
that is $\mathcal{F}_S$-measurably differentiable at $F[S \mid \mathcal{F}_S]$. 
Then $\underline{f} \circ F$  has a stopping derivative at $S$
and 
\begin{equation}
\label{equation:chain}
(\underline{f} \circ F)'\, [\mathcal{F}_S] ~ = ~ 
\underline{f}'(F[S \mid \mathcal{F}_S])\, F'[\mathcal{F}_S] \enspace .
\end{equation}
\end{proposition}

\begin{proof}
Consider any $\mathcal{F}_S$-measurable $\underline{\varepsilon} > 0$. Define:
\[
\RVeps_1 := \frac{\RVeps}{3 \,(|F'[\mathcal{F}_S] | \vee 1)}
\qquad
\RVeps_2 = \frac{\RVeps}{3\, (|\underline{f}'(F[S \mid \mathcal{F}_S])| \vee 1)}
 \enspace .
\]

By the differentiability of $\underline{f}$ at $F[S \mid \mathcal{F}_S]$, let $\underline{\delta} > 0$ be an $\mathcal{F}_S$-measurable
random variable such that, for $\mathcal{F}_S$-almost-all $\omega$ and all $x$ with $ |x - F[S \mid \mathcal{F}_S](\omega)| \leq \underline{\delta}(\omega)$,
\begin{equation}
\label{equation:chain-z}
\left| \underline{f}_\omega'(F[S \mid \mathcal{F}_S](\omega)) - \frac{\underline{f}_\omega (x) - \underline{f}_\omega (F[S \mid \mathcal{F}_S](\omega))}{x - F[S \mid \mathcal{F}_S](\omega)} \right| 
~ \leq ~ \RVeps_1(\omega) \enspace ,
\end{equation}
where, for convenience below, we take the value of the left-hand side to be $0$ in the case that $x = F[S \mid \mathcal{F}_S](\omega)$ and so the numerator and denominator of the fraction are  both $0$.

Using the stopping derivative property of $F'[\mathcal{F}_S] $, let $U' > S$ be such that, 
for any finite $T$ with $S \leq  T \leq U'$,
\begin{equation}
\label{equation:chain-a}
\left(F'[\mathcal{F}_S]  - \RVeps_2\right) \Exp[T-S \mid \mathcal{F}_S] ~ \leq ~ 
F[T \mid \mathcal{F}_S] - F[S \mid \mathcal{F}_S]
~ \leq ~\left(F'[\mathcal{F}_S]  + \RVeps_2\right) \Exp[T-S \mid \mathcal{F}_S]
 \enspace ,
\end{equation}
Define
\[
U~ := ~ U' \wedge \left(S + \frac{\underline{\delta}}{|F'[\mathcal{F}_S]| + \RVeps_2}\right) \enspace .
\]
Clearly, $U$ is a a stopping time with $S < U \leq U'$. Moreover, for any $T$ with
$S \leq T \leq U$, 
\[\Exp[T-S \mid \mathcal{F}_S] ~ \leq ~ \frac{\underline{\delta}}{|F'[\mathcal{F}_S]| + \RVeps_2} \enspace .
\]
So it follows from~\eqref{equation:chain-a} that
\[
|F[T \mid \mathcal{F}_S] - F[S \mid \mathcal{F}_S]| \leq \underline{\delta} \enspace .
\]
Thus by~\eqref{equation:chain-z}, we have:
\[
\left| \underline{f}'(F[S \mid \mathcal{F}_S]) - \frac{\underline{f} (F[T \mid \mathcal{F}_S] ) - \underline{f} (F[S \mid \mathcal{F}_S])}{F[T \mid \mathcal{F}_S] - F[S \mid \mathcal{F}_S]} \right| 
~ \leq ~ \RVeps_1 \enspace ,
\]
again understanding the value of the left-hand side as $0$ in the case of $T$ and $\omega$  for which $F[T \mid \mathcal{F}_S](\omega) = F[S \mid \mathcal{F}_S](\omega)$.
Hence, using~\eqref{equation:chain-a}, we obtain
\begin{align*}
 & \left(\underline{f}'(F[S \mid \mathcal{F}_S])\, F'[\mathcal{F}_S] 
- |F'[\mathcal{F}_S]|\,\RVeps_1  - |\underline{f}'(F[S \mid \mathcal{F}_S])|\, \RVeps_2  - \RVeps_1\RVeps_2 \right) \, \Exp[T-S \mid \mathcal{F_S}]
\\
& \qquad  \leq ~\underline{f} (F[T \mid \mathcal{F}_S] ) - \underline{f} (F[S \mid \mathcal{F}_S]) \\
& \qquad \qquad  \leq ~ 
\left(\underline{f}'(F[S \mid \mathcal{F}_S])\, F'[\mathcal{F}_S] 
+ |F'[\mathcal{F}_S]|\,\RVeps_1 + |\underline{f}'(F[S \mid \mathcal{F}_S])|\, \RVeps_2   + \RVeps_1\RVeps_2\right) \, \Exp[T-S \mid \mathcal{F_S}]  \enspace .
\end{align*}
By the definitions of $\RVeps_1$ and $\RVeps_2$, we have:
\[
|F'[\mathcal{F}_S]|\,\RVeps_1 + |\underline{f}'(F[S \mid \mathcal{F}_S])|\, \RVeps_2   + \RVeps_1\RVeps_2 ~ \leq ~ \frac{7}{9}\,\RVeps  \enspace .
\]
It thus follows that:
\begin{align*}
 & \left( \underline{f}'(F[S \mid \mathcal{F}_S])\, F'[\mathcal{F}_S] 
-  \RVeps \right) \, \Exp[T-S \mid \mathcal{F_S}] ~ 
\\ & \qquad
\leq ~\underline{f} (F[T \mid \mathcal{F}_S] ) - \underline{f} (F[S \mid \mathcal{F}_S]) 
~  
\\ & \qquad \qquad
\leq ~ 
\left(\underline{f}'(F[S \mid \mathcal{F}_S])\, F'[\mathcal{F}_S] 
+  \RVeps\right)\, \Exp[T-S \mid \mathcal{F_S}]  \enspace ,
\end{align*}
which establishes the required equality ~\eqref{equation:chain}.
\end{proof}

A  second version of the chain rule applies to stopping derivatives under 
change of time. 
A \emph{time-change} relative to the filtration $(\mathcal{F})$ is an increasing family $(R_s)_{s \geq 0}$
of finite $\mathcal{F}$-stopping times 
such that, for every $\omega \in \Omega$, the trajectory $s \mapsto R_s(\omega)$ is right continuous.
(Of course  $(R)$ is not required to be $(\mathcal{F})$-adapted.) 
Such a time-change determines a \emph{time-changed filtration} $( \mathcal{F}_{R_s})_{s\geq 0}$, which we abbreviate as $( \mathcal{F}_{R})$.
Due to the stopping-time requirement on $(R)$, for any $(\mathcal{F})$-progressive process $(X)$, the
process $(X_R) := (X_{R_s})_{s \geq 0}$ is  $( \mathcal{F}_{R})$-adapted. Instantiating this in the case of
$(X) := (t)_{t \geq 0}$, we obtain that  $(R)$ is itself  an $( \mathcal{F}_{R})$-adapted process, hence 
$( \mathcal{F}_{R})$-progressive by right continuity. 
If $S$ is a finite stopping time relative to $( \mathcal{F}_{R})$
 then $R_S$ is a stopping time relative to $(\mathcal{F})$. 
 Similarly, if
 $F$ is a stopping functional relative to  $(\mathcal{F})$, then 
 \[
 F \circ (R) ~ := ~ S \leq S' ~ \mapsto ~ F[\,R_{S'} \mid \mathcal{F}_{R_S} \,]\enspace ,
 \]
 is a stopping functional relative to $( \mathcal{F}_{R})$. 
 
 
The time-change rule below is the appropriate form of the chain rule for calculating stopping derivatives  at a stopping time $S$ under a \emph{continuous} time-change $(R)$ in the case that
$(R)$ has a \emph{right path-derivative} at $S$. 

\begin{definition}[Right path-derivative]
\label{definition:path-derivative}
An $\mathcal{F}_S$-measurable random variable $D_S$ is a  
 \emph{right path-derivative} for an $\Real$-valued  $\mathcal{F}_S$-progressive 
 process $(X)$ {at a stopping time $S < \infty$}  if, for
 $\mathcal{F}_S$-almost-all $\omega \in \Omega$, the function $t \mapsto X_t(\omega)$ has right derivative $D_{S}(\omega)$ at
 $S(\omega)$. If such a right path-derivative exists, we say that $(X)$ is \emph{right path-differentiable at $S$}.
 
 A progressive $\Real$-valued process $(D)$ is a  \emph{right-path-derivative process} for $(X)$ if, for every stopping 
 time $S < \infty$, it holds that
 $D_S$ is a right path-derivative for $(X)$ at $S$. If such a right-path-derivative process exists, we say that $(X)$ is \emph{right path-differentiable}.
 \end{definition}

\begin{proposition}[Time-change rule]
\label{proposition:time-change}
Suppose $(R)$ is a continuous time-change
and  $F$ a stopping functional 
relative to $(\mathcal{F})$.
Let $S$ be a finite $(\mathcal{F}_R)$-stopping time such that
$(R)$ has right-path-derivative $D_S$ at $S$.
Suppose further that $F$ has a stopping derivative at $R_S$.
Then the $( \mathcal{F}_R)$-stopping functional
$F \circ (R)$ has stopping derivative at $S$ and
\begin{equation}
\label{equation:time-change}
(F \circ (R))'\, [\mathcal{F}_{R_S}] ~ = ~ 
F'[\mathcal{F}_{R_S}] \cdot D_S \enspace .
\end{equation}
\end{proposition}

\begin{proof}
Consider any $\mathcal{F}_{R_S}$-measurable $\underline{\varepsilon} > 0$. Define:
\[
\RVeps_1 := \frac{\RVeps}{3 \,(D_S  \vee 1)}
\qquad
\RVeps_2 = \frac{\RVeps}{3\, (|F'[\mathcal{F}_{R_S}] |\vee 1)} \enspace .
\]
(Note that $D_S \geq 0$, by the increasing property of $(R)$.)
Clearly $\RVeps_1 , \RVeps_2 > 0$.

Using the stopping derivative property of $F'[\mathcal{F}_{R_S}] $, let $U'_1 > R_S$ be 
an $(\mathcal{F})$-stopping time such that, 
for any finite $(\mathcal{F})$-stopping time $T$ with $R_S \leq  T \leq U'_1$,
\begin{equation}
\label{equation:time-change-a}
\left(F'[\mathcal{F}_{R_S}]  - \RVeps_1\right) \Exp[T-R_S \mid \mathcal{F}_{R_S}] ~ \leq ~ 
F[T \mid \mathcal{F}_{R_S}] - F[R_S \mid \mathcal{F}_{R_S}]
~ \leq ~\left(F'[\mathcal{F}_{R_S}]  + \RVeps_1\right) \Exp[T-R_S \mid \mathcal{F}_{R_S}]
 \enspace ,
\end{equation}
Define:
\[
U_1(\omega)~ := ~ \inf\left\{s' > S(\omega) ~ \left| ~  R_{s'}(\omega)  \geq U'_1(\omega) \right. \right\} \enspace .
\]
Because $(R)$ is increasing and  right continuous, $U_1$ is an $(\mathcal{F}_R)$-stopping time. 
Moreover,  because $(R)$ is right continuous and $U'_1 > R_S$, we have $U_1  > S$.

%

Let $W \in \mathcal{F}_{R_S}$ be a probability $1$ set containing only   $\omega$ such that
  $s \mapsto R_s(\omega)$ has right derivative $D_{S}(\omega)$ at
 $R_S(\omega)$ and $\RVeps_2(\omega)  > 0$.
 Define:
 \begin{equation}
 \label{equation:time-change-b}
 U_2(\omega) := \begin{cases}
 \inf\left\{s' > S(\omega) ~ \left| ~  \left|\frac{R_{s'}(\omega) - R_{S}(\omega)}{s' - S(\omega)} - D_{S}(\omega)\right|  \right.  \geq \RVeps_2(\omega) \right\} & \text {if $\omega \in W$,} 
 \\
 \infty & \text{otherwise.} \\
\end{cases}
 \end{equation}
 Then $U_2$ is an $(\mathcal{F}_R)$-stopping time by a first-approach-style argument, using that  $(R)$ is right continuous and increasing. 
 By the right-derivative property of $D_S$, it holds that, for all $\omega \in W$, there exists $u_\omega > S(\omega)$ such that, 
 for all $s'$ with $S(\omega) < s' \leq u_\omega$,
 \[
 \left| \frac{R_{s'}(\omega) - R_{S}(\omega)}{s'-S(\omega)} - D_S(\omega) \right| 
<  \RVeps_2(\omega) \enspace .
\]
 So $U_2(\omega) \geq u_\omega > S(\omega)$, for all $\omega \in W$. That is, $U_2 > S$.
 
%
 
 Define $U := U_1 \wedge U_2$.  Then $U$ is an $(\mathcal{F}_R)$-stopping time with $U > S$.
 Consider any finite $(\mathcal{F}_R)$-stopping time $S'$ with $S \leq S' \leq U$. We show below that
 \begin{align}
 \nonumber
 & 
 (F'[\mathcal{F}_{R_S}] \cdot D_S - \RVeps)\, \Exp[S'-S \mid \mathcal{F}_S] \\
 \nonumber
 & \qquad \leq ~ (F \circ (R))[S' \mid \mathcal{F}_{R_S}] - (F \circ (R))[S \mid \mathcal{F}_{R_S}] \\
 & \qquad \qquad
 \label{equation:time-change-goal}
 \leq ~ (F'[\mathcal{F}_{R_S}] \cdot D_S + \RVeps) \,\Exp[S'-S \mid \mathcal{F}_S] \enspace ,
 \end{align}
 which establishes~\eqref{equation:time-change} as required.
 
By definition:
\begin{align}
\label{equation:time-change-def}
 & (F \circ (R))[S' \mid \mathcal{F}_{R_S}] - (F \circ (R))[S \mid \mathcal{F}_{R_S}]
 ~  = ~ F[R_{S'} \mid \mathcal{F}_{R_S}] - F[R_{S} \mid \mathcal{F}_{R_S}] \enspace .
\end{align}
For the $(\mathcal{F})$-stopping time $R_{S'}$, it holds that $R_S \leq R_{S'}$, because $(R)$ is increasing.
We show that also $R_{S'} \leq U'_1$. Because  $(R)$ is increasing and $R_S < U'_1$, we have
\[U_1(\omega) ~ = ~ \sup 
\left\{s' \geq S(\omega) ~ \left| ~  R_{s'}(\omega)  < U'_1(\omega) \right. \right\} \enspace ,
\]
for almost all $\omega$.
Since $(R)$ is continuous, $s' \leq U_1(\omega)$ implies 
$R_{s'}(\omega)  \leq U'_1(\omega)$. Hence, because $S' \leq U_1$, it follows
that $R_{S'} \leq U'_1$.

Since $R_S \leq R_{S'} \leq U'_1$, we have, by~\eqref{equation:time-change-a}, that 
\begin{align}
\nonumber
& \left(F'[\mathcal{F}_{R_S}]  - \RVeps_1\right)\, \Exp[R_{S'} -R_S \mid \mathcal{F}_{R_S}] ~
 \\ & \qquad 
 \nonumber
\leq ~ 
F[R_{S'} \mid \mathcal{F}_{R_S}] - F[R_S \mid \mathcal{F}_{R_S}] ~
 \\ & \qquad \qquad 
 \label{equation:time-change-c}
\leq ~\left(F'[\mathcal{F}_{R_S}]  + \RVeps_1\right)\, \Exp[R_{S'} -R_S \mid \mathcal{F}_{R_S}]
\enspace .
\end{align}
Moreover, since
\[
 \Exp[R_{S'} -R_S \mid \mathcal{F}_{R_S}] ~ = ~ 
 \Exp \left[\left.\left(\frac{R_{S'} -R_S}{S'-S}\right)(S'-S)  \, \right| \, \mathcal{F}_{R_S}\right] \enspace ,
\]
$S' \leq U_2$ and $(R)$ is continuous,
it follows from~\eqref{equation:time-change-b} that
\begin{equation}
\label{equation:time-change-d}
(D_S - \RVeps_2) \, \Exp[S' -S \mid \mathcal{F}_{R_S}]
~ \leq ~
\Exp[R_{S'} -R_S \mid \mathcal{F}_{R_S}]
~ \leq ~
(D_S + \RVeps_2) \, \Exp[S' -S \mid \mathcal{F}_{R_S}] \enspace .
\end{equation}
Combining~\eqref{equation:time-change-c} and~\eqref{equation:time-change-d}, we get
\begin{align}
\nonumber
& \left(F'[\mathcal{F}_{R_S}] \cdot D_S  -  D_S\, \RVeps_1 - |F'[\mathcal{F}_{R_S}]|\, \RVeps_2 - \RVeps_1\,\RVeps_2 \right) \, \Exp[S'- S \mid \mathcal{F}_{R_S}] ~
 \\ & \qquad 
 \nonumber
\leq ~ 
F[R_{S'} \mid \mathcal{F}_{R_S}] - F[R_S \mid \mathcal{F}_{R_S}] ~
 \\ & \qquad \qquad 
 \label{equation:time-change-e}
\leq ~ \left(F'[\mathcal{F}_{R_S}] \cdot D_S  + D_S\, \RVeps_1 + |F'[\mathcal{F}_{R_S}]|\, \RVeps_2 + \RVeps_1\,\RVeps_2 \right) \, \Exp[S'- S \mid \mathcal{F}_{R_S}]
\enspace .
\end{align}
By the definitions of $\RVeps_1$ and $\RVeps_2$, we have:
\[
D_S\, \RVeps_1 + |F'[\mathcal{F}_{R_S}]|\, \RVeps_2 + \RVeps_1\,\RVeps_2 ~ \leq ~ 
\frac{7}{9}\,\RVeps \enspace .
\]
Thus~\eqref{equation:time-change-goal} follows from~\eqref{equation:time-change-def} and~\eqref{equation:time-change-e}.
\end{proof}

\noindent
It would be nice to obtain a generalisation of Proposition~\ref{proposition:time-change} in which $(R)$  is only assumed to be right-continuous. However, the current framework of definitions does not appear to support this additional generality.

\section{Drift}
\label{section:drift}

In this section, we investigate properties of \emph{drift} defined as a stopping derivative as in Definition~\ref{definition:drift}.
For convenience, we use the derivative notation $\Exp'[(X) | \mathcal{F}_S]$ for the drift of $(X)$ at $S$, and
$(\Exp'[(X) | \mathcal{F}_t])_{t \geq 0}$ (more concisely $(\Exp'[(X) | \mathcal{F}])$) for 
a drift process for $(X)$, which, as in Definition~\ref{definition:stopping-derivative},  is assumed to be right-continuous.

\begin{proposition}
\label{proposition:drift-sum}
If $(X)$ and $(Y)$ have drift at $S$ then so does $(X+Y)$ and
\[
\Exp'[(X+Y) | \mathcal{F}_S] ~ = ~ \Exp'[(X) | \mathcal{F}_S] + \Exp'[(Y) | \mathcal{F}_S] \enspace .
\]
\end{proposition}
\begin{proof} 
Since $\Exp[(X +Y)_{T-}  \mid \mathcal{F}_S] =  \Exp[X_{T-} + Y_{T-} \mid \mathcal{F}_S]
= \Exp[X_{T-} \mid \mathcal{F}_S] + \Exp[Y_{T-} \mid \mathcal{F}_S]$, the result is a consequence of 
the linearity of stopping derivatives.
\end{proof}

We shall consider two orthogonal classes of processes that have drift at a stopping time $S$. The first is given by processes whose paths are right differentiable at $S$ (recall Definition~\ref{definition:path-derivative}).

 \begin{proposition}
 \label{proposition:path-differentiable}
 If $D_S$ is a  
 {right path-derivative} for a cadlag process $(X)$  at $S$, then $D_S$ is a drift for $(X)$ at $S$. 
 \end{proposition}
 
\begin{proof}
Let $D_S$ be a  
 {right path derivative} for $(X)$  at $S$. Consider any $\mathcal{F}_S$-measurable $\underline{\varepsilon} > 0$. 
 Let $W \in \mathcal{F}_S$ be a probability $1$ set containing only   $\omega$ such that
  $t \mapsto X_t(\omega)$ has right derivative $D_{S}(\omega)$ at
 $S(\omega)$ and $\underline{\varepsilon}(\omega)  > 0$.
 Define:
 \[
 U(\omega) := \begin{cases}
 \inf\left\{t > S(\omega) ~ \left| ~  \left|\frac{X_{t \pm}(\omega) - X_{S}(\omega)}{t - S(\omega)} - D_{S}(\omega)\right|  \right.  \geq \underline{\varepsilon}(\omega) \right\} & \text {if $\omega \in W$,} 
 \\
 \infty & \text{otherwise,} \\
\end{cases}
 \]
 using the notation introduced in~\eqref{equation:plus-minus}.
 Then $U$ is a stopping time by a first-approach argument, using the cadlag property of $(X)$. We show that $S(\omega) < U(\omega)$,  for all $\omega \in W$.
 By the right-derivative property of $D_S$, it holds that, for all $\omega \in W$, there exists $u_\omega > S(\omega)$ such that, 
 for all $t$ with $S(\omega) < t < u_\omega$,
 \[\left| \frac{X_t(\omega) - X_{S(\omega)}(\omega)}{t-S(\omega)} - D_S(\omega) \right| <  \frac{\underline{\varepsilon}(\omega)}{2} \enspace ,\]
 hence also
 \[\left| \frac{X_{t\!-}(\omega) - X_{S(\omega)}(\omega)}{t-S(\omega)} - D_S(\omega) \right| \leq  \frac{\underline{\varepsilon}(\omega)}{2}  < \underline{\varepsilon}(\omega) \enspace .\]
 So $U(\omega) \geq u_\omega > S(\omega)$, for all $\omega \in W$. That is, $U > S$.
 
 Let $T$ be a finite stopping time with $S \leq T \leq U$. Then it is immediate from the definition of $U$ that, for all $\omega \in W$,
 \[
(D_S(\omega) - \underline{\varepsilon}(\omega)\,) (T(\omega)-S(\omega)) ~ \leq ~ (X_{T-} - X_{S})(\omega) ~ \leq ~ 
(D_S(\omega) + \underline{\varepsilon}(\omega)\,) (T(\omega)-S(\omega)) \enspace .
 \]
 It follows that
 \[
 \Exp[(D_S - \underline{\varepsilon}\,) (T-S) \mid \mathcal{F}_S] ~ \leq ~ 
 \Exp[X_{T\!-} - X_{S}  \mid \mathcal{F}_S] ~ \leq ~ 
 \Exp[(D_S + \underline{\varepsilon}\,) (T-S) \mid \mathcal{F}_S] \enspace .
 \]
 Using the $\mathcal{F}_S$-measurability of $D_S$, $\,\underline{\varepsilon}$ and $X_S$, we obtain the desired inequalities:
 \[
 (D_S - \underline{\varepsilon}\,) \Exp[ T-S \mid \mathcal{F}_S] ~ \leq ~ 
 \Exp[X_{T\!-}   \mid \mathcal{F}_S] - X_{S} ~ \leq ~ 
(D_S + \underline{\varepsilon}\,)  \Exp[ T-S \mid \mathcal{F}_S] \enspace .
 \]
\end{proof}
 
 The second class of processes with drift   we consider are those for which the drift is $0$. 

\begin{definition}[Zero drift process]  We say that a process $(X)$ has \emph{zero drift} if it has drift $0$  at all finite stopping times.
\end{definition}

 \begin{theorem}
 \label{theorem:martingale-driftless}
 Every continuous local martingale has zero drift.
 \end{theorem}
 \begin{proof}
 Suppose $(X)$ is a continuous local martingale. 
The assumption of  continuity allows us to 
replace $X_{T-}$ with $X_T$ in the definition of drift. 
By the local martingale property, 
there exists 
an increasing sequence $(V_n)$  of bounded stopping times such that 
$(V_n) \uparrow \infty$ and every stopped process  $X^{V_n}$  is a martingale. Let $b_n$ be an upper bound for $V_n$.


Let $S$ be any finite stopping time.
Define 
\[
 U := \min\{V_n \mid n \in \Nat, \, V_n > S\}  \enspace .
\]
Clearly $U$ is a finite stopping time. Consider any stopping time $T$ with $S \leq T \leq U$. 
Let $P_n$ be the event that $n$ is the smallest natural number such that $V_n > S$. Then $\{P_n \mid n \in \Nat\}$ is an $\mathcal{F}_S$-partition of $\Omega$, and $P_n$ implies $S \leq T \leq U = V_n \leq b_n$. So:
\begin{align*}
\Exp[X_{T} |  \mathcal{F}_S] & = \Exp[X^U_{T} |  \mathcal{F}_S]
\\
& =  E\left[ \left. \left( \sum_n \,  \One_{P_n}  \cdot X^{V_n}_{T \wedge b_n}\right)  \right|  \mathcal{F}_{S}\right]  \\
& = \sum_n \,  \One_{P_n} \cdot \Exp[X^{V_n}_{T \wedge b_n}  |  \mathcal{F}_{S \wedge b_n}] 
& & \text{by~\eqref{equation:countable-sum}}
\\
& = \sum_n \, \One_{P_n} \cdot X_{S \wedge b_n}  & & \text{by the optional stopping theorem}
\\
& = \sum_n \, \One_{P_n} \cdot X_{S}
\\ & = X_{S} \enspace ,
\end{align*}
where the optional stopping theorem applies, since the stopped process $X^{V_n}_{t}$  is a martingale and
the stopping time $T \wedge b_n$ is bounded.
Since  $\Exp[X_{T} | \mathcal{F}_S] = X_S$, for every $T$ with $S \leq T \leq U$, the drift at $S$  is $0$.
 \end{proof}
 \noindent
 Note that the use of~\eqref{equation:countable-sum} in the above proof means that it is crucial that the definition of drift is defined using generalised conditional expectation.
 
 
 A wide class of processes can be decomposed as a sum of a zero-drift process and a process with right path-derivative. This class of processes includes those that have a \emph{path-integrable} drift process.
 \begin{definition}[Path-integrable process]
 \label{definition:path-integrable}
 We say that a process $(B)$ is \emph{path-integrable} if it is right continuous and,
 for every $\omega \in \Omega$ and $t > 0$,
\begin{equation}
\label{eq:pi-i}
 \int_0^t |B_s(\omega)| \, ds < \infty \enspace .
\end{equation}
If $(B)$ is path-integrable, define its \emph{path-integral process} $(\int \! B)$ as
\begin{equation}
\label{eq:pi-ii}
t , \omega \mapsto \int_0^t B_s(\omega) \, ds \enspace .
\end{equation}
 \end{definition}
 \noindent
 The assumption of right continuity, in the definition of path integrability, could be weakened to progressivity, and~\eqref{eq:pi-ii} would still define $(\int \! B)$ as a progressive process. 
 However,  in the present article, the property of right continuity will always be needed in combination with path integrability, so we include it in the definition. 
 
  \begin{proposition}
 \label{proposition:integrable-gives-drift}
 If $(B)$ is {path-integrable} then its path-integral process $(\int \! B)$ is right-path-differentiable with $(B)$
 as its right-path-derivative process hence also as its drift process.
 \end{proposition}
 \begin{proof}
Because $(B)$ is right continuous, the right-differentiable version of the fundamental theorem of calculus gives us that,
for every $\omega \in \Omega$, the path $t \mapsto (\int \! B)_t (\omega)$ is
continuous and right differentiable, with right derivative
$B_s(\omega)$ at $s$. Accordingly, it follows from
Proposition~\ref{proposition:path-differentiable} that $(B)$ is a drift
process for the process $(\int \! B)$. 
\end{proof}
\begin{corollary}
\label{corollary:kill-drift}
Suppose $(X)$ has a path-integrable drift process $(B)$. 
Then $(X-\int\! B)$ is a zero-drift process. 
\end{corollary}
\begin{proof}
 The process $(X-\int\! B)$  has drift $0$, by
Propositions~\ref{proposition:integrable-gives-drift} and~\ref{proposition:drift-sum}.
\end{proof}

\section{Variance rate}
\label{section:diffusion}

In this section, we investigate properties of the  \emph{variance rate} defined as a stopping derivative as in Definition~\ref{definition:variance-rate}.
We use the derivative notation $\Var'[(X) | \mathcal{F}_S]$ for the variance rate of $(X)$ at $S$ (if it exists), and
$(\Var'[(X) | \mathcal{F}_t])_{t \geq 0}$ (more concisely $(\Var'[(X) | \mathcal{F}])$) for 
the  variance-rate process for $(X)$, which is required to be right-continuous.

The variance rate, if it exists, is always non-negative. Although this is  a very basic observation, we give the proof to  
highlight the convenience of using a random $\RVeps$. 
\begin{proposition}
\label{proposition:non-negative}
If the variance rate of $(X)$ at $S$ exists then $\Var'[(X) | \mathcal{F}_S] \geq 0$.
\end{proposition}
\begin{proof} Suppose for contradiction that  $\Var'[(X) | \mathcal{F}_S] < 0$ on a positive measure $W \subseteq \mathcal{F}_S$. Define:
\[
\RVeps(\omega) ~ := ~ \begin{cases} -\frac{1}{2}\Var'[(X) | \mathcal{F}_S] (\omega) & \text{if $\omega \in W$} \\
1 & \text{otherwise.}
\end{cases}
\]
Clearly $\RVeps$ is $\mathcal{F}_S$-measurable and $\RVeps > 0$. By the definition of variance rate, there exists a finite stopping time $U > S$ such that, for all finite $T$ with $S \leq T \leq U$,
\[
\Var[X_{T-} | \mathcal{F}_S]  ~ \leq ~ 
 (\Var'[(X) | \mathcal{F}_S]  + \underline{\varepsilon}\,) \, \Exp[T-S  \mid \mathcal{F}_S] \enspace .
 \]
 In the case $T := U$, we have $\Var[X_{T-} | \mathcal{F}_S] \geq 0$. Also,  by the choice  of $\RVeps$, for all $\omega \in W$, we have
 $ ((\Var'[(X) | \mathcal{F}_S]  + \underline{\varepsilon}\,) \, \Exp[T-S  \mid \mathcal{F}_S])(\omega) < 0$.
 This contradicts the inequality above.
\end{proof}

In practice it is often convenient to work with a characterisation of the variance rate as the
\emph{relative-second-moment rate} as defined below.


\begin{definition}[Relative-second-moment rate]
The \emph{relative-second-moment rate} of a cadlag process $(X)$ is a stopping derivative for the \emph{relative second moment} function
\begin{equation}
\label{equation:rsm-rate}
\Exp[(X_{T-}  - X_S)^2 \,| \mathcal{F}_S]\, \enspace .
\end{equation}
\end{definition}


\begin{proposition}
\label{proposition:rsm-rate}
Suppose a cadlag process $(X)$ has drift at a stopping time $S < \infty$. Then the following are equivalent for
an $\mathcal{F}_S$-measurable random variable $C_S$.
\begin{enumerate}
\item \label{vrsm:a} $C_S$ is the variance rate for $(X)$ at $S$.
\item \label{vrsm:b} $C_S$ is the relative-second-moment rate for $(X)$ at $S$.
\end{enumerate}
\end{proposition}
\begin{proof}
Let $D_S$ be the drift of $(X)$ at $S$. 

For the conditional variance function, we have:
\[
\Var[X_{T-} \mid \mathcal{F}_S]  ~ = ~ 
\Exp[(X_{T-})^2  \mid \mathcal{F}_S] - 
(\Exp[X_{T-}  \mid \mathcal{F}_S])^2 \enspace .
\]
By the product rule (Proposition~\ref{proposition:product}), the stopping derivative of $(\Exp[X_{T-}  | \mathcal{F}_S])^2$ at $S$  is $2 D_S X_S$.
It then follows from the equation above and linearity (Proposition~\ref{proposition:linear}), that $(X)$ has
variance rate at $S$ if and only if $\Exp[(X_{T-})^2  | \mathcal{F}_S]$ has stopping derivative at $S$, in which case
\begin{equation}
\label{equation:var-square-general}
\Var'[(X)   \mid  \mathcal{F}_S] ~ = ~ 
\Exp'[(X^2)  \mid \mathcal{F}_S] - 2  D_S X_S \enspace .
\end{equation}

Similarly, for the relative second moment function, we have:
\[
\Exp[(X_{T-}  - X_S)^2 \mid \mathcal{F}_S] ~ = ~ 
\Exp[(X_{T-})^2  \mid \mathcal{F}_S] - 2 X_S \Exp[X_{T-}  \mid \mathcal{F}_S]  + (X_S)^2 \enspace .
\]
So, by the linearity and product rules for stopping derivatives,
the relative-second-moment rate exists if and only if $\Exp[(X_{T-})^2  \mid \mathcal{F}_S]$ has a stopping derivative at $S$, in which case
\[
\Exp'[((X  - X_S)^2) \mid \mathcal{F}_S] ~ = ~ 
\Exp'[(X^2) \mid \mathcal{F}_S] - 2 D_S X_S   \enspace .
\]

It follows that each of $\Var'[(X)   \mid  \mathcal{F}_S]$ and $\Exp'[((X  - X_S)^2) \mid \mathcal{F}_S]$ exists if and only if the other does, 
and that when they do exist they are (almost surely) equal.
\end{proof}

As a special case of~\eqref{equation:var-square-general}  in the argument above, we note that if $(X)$ has zero drift at $S$ then
\begin{equation}
\label{equation:var-square}
\Var'[(X)   \mid  \mathcal{F}_S] ~ = ~ 
\Exp'[(X^2)  \mid \mathcal{F}_S] \enspace .
\end{equation}

\begin{proposition}[Brownian motion]
\label{proposition:brownian}
Standard Brownian motion $(\Brown_t)_t$ has drift process $(0)$ and variance-rate process $(1)$.
\end{proposition}
\begin{proof}
Brownian motion  $(\Brown)$  and the process $(\Brown_t^2 - t)_{t \geq 0}$ are continuous martingales, hence both have drift process $(0)$ by 
Theorem~\ref{theorem:martingale-driftless}. Since the process $(t)$ is obviously right path differentiable with right-path-derivative process $(1)$,
the process $(\Brown^2)$ has drift process $(1)$, by Proposition~\ref{proposition:drift-sum}. 
Then, by~\eqref{equation:var-square} above,  $(\Var'[(\Brown)   \mid  \mathcal{F}]) = (\Exp'[(\Brown^2)  \mid \mathcal{F}]) = (1)$.
%
\end{proof}
\noindent
We shall show in Section~\ref{section:fundamental} that the properties of having drift process $(0)$ and variance-rate process $(1)$ together characterise Brownian motion amongst continuous processes.

We now turn to our formulation of It\^o's formula,
a first version of which  was given as \eqref{equation:ito} in Section~\ref{section:intro}. We extend the statement given there in two directions. Firstly, we generalise the time $s$ to a finite stopping time $S$. Secondly, we generalise the (deterministic) function $f$ to a {random function}
$\underline{f}\,$, in the sense of Definition~\ref{definition:random-function}. We shall need, in particular, the notion
of measurable continuous differentiability for random functions from Definition~\ref{definition:measurable-continuous-differentiability}.

If $S$ is a stopping time and $\underline{f}$ is an $\mathcal{F}_S$-jointly-measurable random function, then we write
$(\underline{f}(X_t))_{t \geq S}$ for the process 
\begin{equation}
\label{equation:restrict-S}
t, \omega \mapsto \begin{cases} \underline{f}_\omega(X_t(\omega)) & \text{if $t \geq S(\omega)$} \\
      0 & \text{otherwise,}
      \end{cases}
\end{equation}
or simply $(\underline{f}(X))$ when $S$ is clear from the context. 
Note that if the random function $\underline{f}$ is continuous and $(X)$ is cadlag then $(\underline{f}(X))$ is also cadlag.

\begin{theorem}[The It\^o formula]
\label{theorem:ito-formula}
Suppose that  $(X)$ is a cadlag process with drift  and variance rate at a finite stopping time $S$.
Suppose that $\underline{f}$ is an  $\mathcal{F}_S$-jointly-measurable random function such that
$\underline{f}_\omega(X_t(\omega))$ is defined for all $t \geq S(\omega)$.
Suppose further that  $\underline{f}$ is twice
$\mathcal{F}_S$-measurably continuously differentiable at $X_S$.
Then the process $(\underline{f}(X))$ has both drift 
and variance rate at $S$, with these given by:
\begin{align}
\label{equation:ito-exp}
\Exp'[ (\underline{f}(X)) \mid \mathcal{F}_S] ~ = & ~  \underline{f}'(X_S)\,\Exp'[(X) \mid \mathcal{F}_S] +
 \frac{\underline{f}''(X_S)}{2}\,\Var'[(X) \mid \mathcal{F}_S] \enspace , \\
 \label{equation:ito-var}
\Var'[ (\underline{f}(X)) \mid \mathcal{F}_S] ~ = & ~ (\underline{f}'(X_S))^2 \,\Var'[(X) \mid \mathcal{F}_S]\enspace .
\end{align}
\end{theorem}
\begin{proof} Let $B_S$ and $A_S$ be the drift and variance rate of $(X)$ at $S$ respectively.
We first show that $(\underline{f}(X))$ has drift at $S$ satisfying~\eqref{equation:ito-exp}.
Let $\RVeps > 0$ be $\mathcal{F}_S$-measurable. Define:
\[
\RVeps_1 := \frac{\RVeps}{4 \,(|\underline{f}'(X_S)| \vee 1)}
\qquad
\RVeps_2 = \frac{\RVeps}{2\, (|\underline{f}''(X_S)| \vee 1)}
\qquad
\RVeps_3 = \frac{\RVeps}{2\, (A_S \vee 1)} \enspace .
\]
Then $\RVeps_1, \RVeps_2, \RVeps_3$ are $\mathcal{F}_S$-measurable and
clearly $0 < \RVeps_1, \RVeps_2, \RVeps_3$. 

Using the drift property of $B_S$, let $U_1 > S$ be such that, for all finite $T$ with $S \leq T \leq U_1$,
\begin{equation}
\label{equation:AB-a}
(B_S - \varepsilon_1 ) \,\Exp[T-S \mid  \mathcal{F}_S ] ~ \leq ~ \Exp[X_{{T \! -}}| \mathcal{F}_S] - X_S ~ \leq ~ (B_S + \varepsilon_1 ) \, \Exp[T-S \mid  \mathcal{F}_S ] \enspace .
\end{equation}
Similarly, using the relative-second-moment-rate characterisation of the variance rate $A_S$,  let $U_2 > S$ be such that, for all finite $T$ with $S \leq T \leq U_2$,
\begin{equation}
\label{equation:AB-b}
(A_S - \varepsilon_2 ) \, \Exp[T-S \mid  \mathcal{F}_S ] ~ \leq ~ \Exp[(X_{{T \! -}}-X_S)^2 | \mathcal{F}_S]  ~ \leq ~ (A_S + \varepsilon_2 )\, \Exp[T-S \mid  \mathcal{F}_S ] \enspace .
\end{equation}
Let $E(\omega) \subseteq \Real$ be the largest open interval containing $X_S(\omega)$ and satisfying, for
all $x \in E(\omega)$, it holds that $|\underline{f}_\omega''(x)-\underline{f}_\omega''(X_S(\omega))| < \RVeps_3(\omega)$. Such an open interval exists because 
$\underline{f}_\omega''$ is continuous. Define
\[
U_3(\omega)~  := ~  \inf \{t \geq S(\omega) \mid \underline{f}_\omega''(X_{t\pm}(\omega)) \notin E(\omega)\} 
\enspace .
\]
Because  $\underline{f}''$ is a continuous $\mathcal{F}_S$-jointly-measurable random function
and $(X)$ is cadlag, we have  that: 
$U_3$ is a stopping time (by a first-approach argument), $U_3 > S$ and, for all $t$ with
$S(\omega) \leq t  \leq U_3(\omega)$ and all $\xi \in \Real$ between $X_{S}(\omega)$ and 
$X_{t-}(\omega)$ (it is important that this is $X_{t-}(\omega)$ rather that $X_t(\omega)$), 
\begin{equation}
\label{equation:AB-z}
|\underline{f}_\omega''(\xi) - \underline{f}_\omega''(X_S(\omega))| \leq \RVeps_3 \enspace .
\end{equation}

Define $U := U_1 \wedge U_2 \wedge U_3$, and 
consider any finite $T$ with $S \leq T \leq U$. 
We analyse $\underline{f}(X_{T\!-}) - \underline{f}(X_S)$, using a second-order Taylor expansion.
\begin{align*}
& \underline{f}_\omega (X_{T\!-}(\omega)) -  \underline{f}_\omega (X_S(\omega))  ~ = 
\\
& \qquad  
  \underline{f}_\omega'(X_{S}(\omega))(X_{T\!-}(\omega) \!- \!X_{S}(\omega)) + 
 \underline{f}_\omega''(X_{S}(\omega))\left(\frac{(X_{T\!-}(\omega) \!- \!X_{S}(\omega))^2}{2}\right) + R_\omega \enspace,
\end{align*}
where the error term $R_\omega$ is given by
\[
R_\omega ~ := ~ (\underline{f}''(\xi_\omega) - \underline{f}''(X_{S}(\omega)))\left(\frac{(X_{T\!-}(\omega) \!- \!X_{S}(\omega))^2}{2}\right) \enspace ,
\]
for some $\xi_\omega$ between $X_{S}(\omega)$ and $X_{T\!-}(\omega)$.
Since  $S(\omega) \leq T(\omega) \leq U_3(\omega)$, it follows from~\eqref{equation:AB-z} that
\[
|R_\omega| ~ \leq ~ \RVeps_3 \left(\frac{(X_{T\!-}(\omega) \!- \!X_{S}(\omega))^2}{2}\right) \enspace .
\]
Accordingly, we have 
\begin{align}
\nonumber
& \underline{f}_\omega'(X_{S}(\omega))(X_{T\!-}(\omega) \!- \!X_{S}(\omega)) + 
 (\underline{f}_\omega''(X_{S}(\omega)) - \RVeps_3)\left(\frac{(X_{T\!-}(\omega) \!- \!X_{S}(\omega))^2}{2}\right) \\
 \nonumber
 &  \qquad \leq ~ \underline{f}_\omega(X_{T\!-}(\omega)) -  \underline{f}_\omega(X_S(\omega)) \\
 \label{equation:AB-y}
 & \qquad \qquad \leq ~
 \underline{f}_\omega'(X_{S}(\omega))(X_{T\!-}(\omega) \!- \!X_{S}(\omega)) + 
 (\underline{f}_\omega''(X_{S}(\omega)) + \RVeps_3)\left(\frac{(X_{T\!-}(\omega) \!- \!X_{S}(\omega))^2}{2}\right)
\end{align}
It follows that:
\begin{align}
\nonumber
& \underline{f}'(X_{S})\Exp[X_{T\!-} \!- \!X_{S} \mid \mathcal{F}_S] + 
 (\underline{f}''(X_{S}) - \RVeps_3)\left(\frac{\Exp[(X_{T\!-}\!- \!X_{S})^2 \mid \mathcal{F}_S]}{2}\right) \\
  \nonumber
 & \qquad \leq ~ \Exp[\underline{f}(X_{T\!-}) | \mathcal{F}_S] - \underline{f}(X_S)\\
 \label{equation:AB-c}
 & \qquad \qquad \leq ~
 \underline{f}'(X_{S})\Exp[X_{T\!-} \!- \!X_{S} \mid \mathcal{F}_S] + 
 (\underline{f}''(X_{S}) + \RVeps_3)\left(\frac{\Exp[(X_{T\!-}\!- \!X_{S})^2 \mid \mathcal{F}_S]}{2}\right) \enspace .
\end{align}
By~\eqref{equation:AB-a}, we have
\begin{align}
\nonumber
& (\underline{f}'(X_{S})B_S - |\underline{f}'(X_{S})|\,\RVeps_1 ) \,\Exp[T-S \mid  \mathcal{F}_S ] \\
\nonumber
& \qquad \leq ~ \underline{f}'(X_{S})\Exp[X_{T\!-} \!- \!X_{S} \mid \mathcal{F}_S] \\
\label{equation:AB-d}
& \qquad \qquad  \leq ~ (\underline{f}'(X_{S})B_S + |\underline{f}'(X_{S})|\,\RVeps_1 ) \, \Exp[T-S \mid  \mathcal{F}_S ] \enspace .
\end{align}
Via a slightly more involved argument using~\eqref{equation:AB-b}, we also have the two inequalities:
\begin{align}
\nonumber
& \left(\frac{\underline{f}''(X_{S})}{2} A_S - \frac{|\underline{f}''(X_{S})|\,\RVeps_2 \!+\! A_S\, \RVeps_3 \!+\! \RVeps_2\RVeps_3}{2}\right) \!\Exp[T\!-\!S | \mathcal{F}_S]
\leq  (f''(X_{S}) - \RVeps_3)\!\left(\frac{\Exp[(X_{T\!-}\!- \!X_{S})^2 | \mathcal{F}_S]}{2}\right) \enspace,\\
 \label{equation:AB-e}
 & (\underline{f}''(X_{S}) + \RVeps_3)\!\left(\frac{\Exp[(X_{T\!-}\!- \!X_{S})^2 | \mathcal{F}_S]}{2}\right)
 \leq 
 \left(\frac{\underline{f}''(X_{S})}{2} A_S + \frac{|\underline{f}''(X_{S})|\,\RVeps_2 \!+ \!A_S \,\RVeps_3 \!+\! \RVeps_2\RVeps_3}{2}\right) \! \Exp[T\!-\!S |\mathcal{F}_S]
  \enspace .
\end{align}
Combining~\eqref{equation:AB-c},~\eqref{equation:AB-d} and~\eqref{equation:AB-e}, we get:
\begin{align}
\nonumber
& \left(\underline{f}'(X_{S})B_S + \frac{\underline{f}''(X_{S})}{2} A_S -
|\underline{f}'(X_{S})|\,\RVeps_1 - \frac{|\underline{f}''(X_{S})|\,\RVeps_2\!+ \!A_S\, \RVeps_3 \!+ \!\RVeps_2\RVeps_3}{2}\right) \Exp[T\!-\!S | \mathcal{F}_S]
\\
\nonumber
 & \quad \leq  \Exp[\underline{f}(X_{T\!-}) | \mathcal{F}_S] - \underline{f}(X_S) \\
 \label{equation:AB-f}
 & \quad \quad \leq 
 \left(\underline{f}'(X_{S})B_S + \frac{\underline{f}''(X_{S})}{2} A_S +
|f'(X_{S})|\,\RVeps_1 + \frac{|f''(X_{S})|\,\RVeps_2\! + \!A_S \,\RVeps_3 \!+ \!\RVeps_2\RVeps_3}{2}\right)
\Exp[T\!-\!S | \mathcal{F}_S]
  \enspace .
\end{align}
By the definitions of $\RVeps_1, \RVeps_2, \RVeps_3$, we have:
\[
|f'(X_{S})|\,\RVeps_1 + \frac{|f''(X_{S})|\,\RVeps_2\! + \!A_S \,\RVeps_3 \!+ \!\RVeps_2\RVeps_3}{2} ~ \leq~  \frac{7}{8}\,\RVeps \enspace .
\]
It thus follows from~\eqref{equation:AB-f} that:
\begin{align*}
& \left(\underline{f}'(X_{S})B_S + \frac{\underline{f}''(X_{S})}{2} A_S - \RVeps\right) \Exp[T\!-\!S | \mathcal{F}_S] \\
& \qquad  \leq  ~ \Exp[\underline{f}(X_{T\!-}) | \mathcal{F}_S] - \underline{f}(X_S) \\
& \qquad \qquad \leq 
 \left(\underline{f}'(X_{S})B_S + \frac{\underline{f}''(X_{S})}{2} A_S + \RVeps\right) \Exp[T\!-\!S | \mathcal{F}_S] \enspace .
 \end{align*}
 This establishes that $(\underline{f}(X_t))_{t \geq S}$ has drift  at $S$ satisfying~\eqref{equation:ito-exp}.
 
 It remains to show that $(\underline{f}(X))$ has variance rate at $S$ satisfying~\eqref{equation:ito-var}.
By the characterisation of variance rate as relative-second-moment rate, $C_S$  is a variance rate for
$(\underline{f}(X))$ at $S$ if and only if it is a drift for the process $(\underline{g}(X))$, where 
$\underline{g}(x) = (\underline{f}(x) - \underline{f}(X_S))^2$.
Since the random function $\underline{g}$ has continuous derivatives 
$\underline{g}'(x) = 2(\underline{f}(x) - \underline{f}(X_S))\underline{f}'(x)$ and
$\underline{g}''(x) = 2(\underline{f}'(x))^2+ 2 (\underline{f}(x) - \underline{f}(X_S))\underline{f}''(x)$, 
we have: 
\[
\Var'[(\underline{f}(X)) | \mathcal{F}_S] ~ = ~ 
\Exp'[(\underline{g}(X)) | \mathcal{F}_S] ~ = ~ \underline{g}'(X_S)B_S + \frac{\underline{g}''(X_S)}{2}A_S ~ = ~ (\underline{f}'(X_S))^2 A_S \enspace ,
\]
using the characterisation of drift by the previously established~\eqref{equation:ito-exp}.
%
%
\end{proof}
\begin{remark}
\label{remark:Ito}
It is worth highlighting the point in the above argument that necessitates using a second-order Taylor expansion 
to obtain the drift formula.  If one attempts to use a first-order Taylor expansion, then the analogues of the inequalities
in~\eqref{equation:AB-y} no longer hold. This is because the nonnegative term $(X_{T\!-}(\omega) \!- \!X_{S}(\omega))^2$
in~\eqref{equation:AB-y} is replaced with $X_{T\!-}(\omega) \!- \!X_{S}(\omega)$, which may be negative or positive depending on $\omega$,
and so there is no  uniform expression for the lower and upper bounds,
applying to all $\omega$, that can be used to continue the argument to obtain an
analogue of \eqref{equation:AB-c}.
\end{remark}

\section{Covariance rate}
\label{section:covariance}

\begin{definition}[Covariance rate]
\label{definition:covariance-rate}
The \emph{covariance rate} between two cadlag processes $(X)$ and $(Y)$ at $S$ is a stopping derivative for the \emph{conditional covariance} function
\[
\Cov[X_{T-}, Y_{T-} \mid \mathcal{F}_S] ~ := ~ \Exp[(X_{T-}  - \Exp[X_{T-} | \mathcal{F}_S])(Y_{T-}  - \Exp[Y_{T-} | \mathcal{F}_S]) \mid \mathcal{F}_S]\, \enspace .
\]
\end{definition}
\noindent
We  use the derivative-style notation $\Cov'[(X),(Y) \mid \mathcal{F}_S]$ for the covariance rate at $S$.

Covariance rate generalises variance rate because $\Var'[(X) \mid \mathcal{F}_S] = \Cov'[(X),(X) \mid \mathcal{F}_S]$.
Given this, the result below and its proof directly generalise  Proposition~\ref{proposition:rsm-rate}.
\begin{proposition}
\label{proposition:cov-rate}
Suppose cadlag processes $(X)$ and $(Y)$ have drift at a stopping time $S < \infty$. Then the following are equivalent for
an $\mathcal{F}_S$-measurable random variable $C_S$.
\begin{enumerate}
\item \label{cvrsm:a} $C_S$ is the covariance rate between $X$ and $Y$ at $S$.
\item \label{vvrsm:b} $C_S$ is a stopping derivative for 
\[\Exp[(X_{T-}  - X_S)(Y_{T-}  - Y_S) \mid \mathcal{F}_S]\, \enspace .
\]
\end{enumerate}
\end{proposition}
\begin{proof}
Let $D_S$ and $E_S$ respectively be drifts of $(X)$ and $(Y)$ at $S$.

For the conditional covariance function, we have the standard covariance equality:
\[
\Cov[X_{T-}, Y_{T-}  \mid \mathcal{F}_S]  ~ = ~ 
\Exp[X_{T-} \, Y_{T-}  \mid \mathcal{F}_S] - 
\Exp[X_{T-}  \mid \mathcal{F}_S] \,\Exp[Y_{T-}  \mid \mathcal{F}_S] \enspace .
\]
By Proposition~\ref{proposition:product}, the stopping derivative of 
$\Exp[X_{T-}  \mid \mathcal{F}_S] \, \Exp[Y_{T-}  \mid \mathcal{F}_S]$
at $S$  is $D_S Y_S + X_S E_S$.
It then follows from the equation above and linearity (Proposition~\ref{proposition:linear}), that $(X)$ and $(Y)$ have
covariance rate at $S$ if and only if $\Exp[X_{T-} \, Y_{T-}  \mid \mathcal{F}_S]$ has stopping derivative at $S$, in which case
\[
\Cov'[(X), (Y)    \mid  \mathcal{F}_S] ~ = ~ 
\Exp'[(XY)  \mid \mathcal{F}_S] -  D_S Y_S - X_S E_S \enspace .
\]

Similarly,  we have:
\[
\Exp[(X_{T-}  - X_S) (Y_{T-} - Y_S)  \mid \mathcal{F}_S] ~ = ~ 
\Exp[X_{T-}\,Y_{T-}  \mid \mathcal{F}_S] - X_S \Exp[Y_{T-}  \mid \mathcal{F}_S]  - Y_S \Exp[X_{T-}  \mid \mathcal{F}_S]  + X_S Y_S\enspace .
\]
So, by the linearity and product rules for stopping derivatives,
 the left-hand side has a stopping derivative if and only if $\Exp[X_{T-}\,Y_{T-}  \mid \mathcal{F}_S]$ does, in which case
\[
\Exp'[((X  - X_S)(Y-Y_S))\mid \mathcal{F}_S] ~ = ~ 
\Exp'[(XY) \mid \mathcal{F}_S] -   D_S Y_S - X_S E_S  \enspace .
\]

It follows that each of $\Cov'[(X), (Y)  \mid  \mathcal{F}_S]$ and $\Exp'[((X  - X_S)(Y-Y_S))\mid \mathcal{F}_S]$ exists if and only if the other does, and that when they do exist they are (almost surely) equal.
\end{proof}

\begin{proposition}
\label{proposition:bilinearity-covariance}
If the pairs $(X), (Z)$ and $(Y),(Z)$ both have covariance rate at $S$ then also  $(X+Y), (Z)$ has covariance rate at $S$ and
\[
\Cov'[(X+Y), (Z)  \mid \mathcal{F}_S] ~ = ~ \Cov'[(X),(Z) \mid \mathcal{F}_S] + \Cov'[(X),(Z) \mid \mathcal{F}_S] \enspace .
\]
\end{proposition}
\begin{proof}
The result follows from the bilinearity of covariance and the linearity of stopping derivatives.
\end{proof}

\begin{proposition}
\label{proposition:variance-sum}
If $(X)$ and $(Y)$ have variance rate  at $S$ then the sum process $(X+Y)$ has variance rate at $S$ if and only if 
the pair $(X),(Y)$ has covariance rate at $S$, in which case  
\[
\Var'[(X+Y) \mid \mathcal{F}_S] ~ = ~ \Var'[(X) \mid \mathcal{F}_S] + 2\, \Cov'[(X),(Y) \mid \mathcal{F}_S] + \Var'[(Y) \mid \mathcal{F}_S]  \enspace .
\]
\end{proposition}
\begin{proof}
A direct consequence of Proposition~\ref{proposition:bilinearity-covariance}, using the general identity $\Var[(Z) \mid \mathcal{F}_S] = \Cov[(Z),(Z) \mid \mathcal{F}_S]$.
\end{proof}

\begin{proposition}
\label{proposition:product-variance}
If $(X)$ and $(Y)$ have drift at $S$ then the product process $(XY)$ has drift at $S$ if and only if 
the pair $(X),(Y)$ has covariance rate at $S$, in which case  
\[
\Exp'[(XY) \mid \mathcal{F}_S] ~ = ~ X_S \,\Exp'[(Y) \mid \mathcal{F}_S] + Y_S \,\Exp'[(X) \mid \mathcal{F}_S] + \Cov'[(X),(Y) \mid \mathcal{F}_S] \enspace .
\]
\end{proposition}

\begin{proof}
The result follows from the standard covariance equality
\[
\Exp[X_{T-} \, Y_{T-}  \mid \mathcal{F}_S]  ~ = ~ 
\Exp[X_{T-}  \mid \mathcal{F}_S]\, \Exp[Y_{T-}  \mid \mathcal{F}_S] + \Cov[X_{T-}, Y_{T-}  \mid \mathcal{F}_S] 
\enspace ,
\]
by linearity and the product rule for stopping derivatives.
\end{proof}

Processes with right path derivatives (see Definition~\ref{definition:path-derivative}) always have covariance rate $0$. 
\begin{proposition}
\label{proposition:rpd-covariance}
If a cadlag process $(X)$ is right-path-differentiable at $S$, then for every cadlag process $(Y)$, it holds that
$\Cov'[(X),(Y) \mid \mathcal{F}_S] = 0$. In particular $\Var'[(X)\mid \mathcal{F}_S] = 0$. 
\end{proposition}
\begin{proof}
Let $D_S$ be a  
 {right path derivative} for $(X)$  at $S$. Consider any $\mathcal{F}_S$-measurable $\underline{\varepsilon} > 0$. 
 Define:
 \[
 \RVeps_1 ~ := ~ 1 \qquad \RVeps_2 ~ := ~ \frac{\RVeps}{|D_S| +1} \enspace .
 \]
 
 Let $W \in \mathcal{F}_S$ be a probability $1$ set containing only   $\omega$ such that
  $t \mapsto X_t(\omega)$ has right derivative $D_{S}(\omega)$ at
 $S(\omega)$.
 Define:
 \[
 U_1(\omega) := \begin{cases}
 \inf\left\{t > S(\omega) ~ \left| ~  \left|\frac{X_{t \pm}(\omega) - X_{S}(\omega)}{t - S(\omega)} - D_{S}(\omega)\right|  \right.  \geq \RVeps_1(\omega) \right\} & \text {if $\omega \in W$,} 
 \\
 \infty & \text{otherwise,} \\
\end{cases}
 \]
 using the notation introduced in~\eqref{equation:plus-minus}.
 As in the proof of Proposition~\ref{proposition:path-differentiable}, 
 $U_1$ is a stopping time and $U_1 > S$.
 Additionally define:
 \[
 U_2(\omega)~  :=  ~
 \inf\left\{t > S(\omega) ~ \mid \; |Y_{t \pm}(\omega) - Y_{S}(\omega)|   \geq \RVeps_2(\omega) \right\} \enspace .
 \]
 Then, because $(Y)$ is right continuous,  $U_2$ is a stopping time with $U_2 > S$. Define $U := U_1 \wedge U_2$.
 
Let $T$ be a finite stopping time with $S \leq T \leq U$. It follows from the definition of $U_1$ that, 
 \[
(D_S - \RVeps_1\,) (T-S) ~ \leq ~ (X_{T-} - X_{S})~ \leq ~ 
(D_S + \RVeps_1\,) (T-S) \enspace ,
 \]
 and from the definition of $U_2$ that $|Y_{T-} - Y_S| \leq \RVeps_2$. Thus we have
\[
-(|D_S| + \RVeps_1\,)\,\RVeps_2\, (T-S) ~ \leq ~ (X_{T-} - X_{S})(Y_{T-} - Y_S) ~ \leq ~ 
(|D_S| + \RVeps_1 \,)\,\RVeps_2\, (T-S) \enspace .
 \]
So, by the definitions of $\RVeps_1$ and $\RVeps_2$,
\[
-\RVeps\, (T-S) ~ \leq ~ (X_{T-} - X_{S})(Y_{T-} - Y_S) ~ \leq ~ 
\RVeps\, (T-S) \enspace .
\]
It follows that
\[
-\RVeps\, \Exp[T-S \mid \mathcal{F}_S] ~ \leq ~ \Exp[(X_{T-} - X_{S})(Y_{T-} - Y_S) \mid \mathcal{F}_S] ~ \leq ~ 
\RVeps\, \Exp[T-S \mid \mathcal{F}_S] \enspace ,
\]
establishing that $\Cov'[(X),(Y) \mid \mathcal{F}_S] = 0$ via Proposition~\ref{proposition:cov-rate}.
\end{proof}

\begin{proposition}
\label{proposition:variance-rate-preserved}
Suppose $(X)$ has drift and variance rate at a finite stopping time $S$, and suppose 
$(Y)$ is right-path differentiable at $S$. Then $(X + Y)$ has variance-rate process at $S$ and
$\Var'[(X + Y) \mid \mathcal{F}_S] = \Var'[(X) \mid \mathcal{F}_S]$.
\end{proposition} 
\begin{proof} We have:
\begin{align*}
& \Var'[(X+Y
) \mid \mathcal{F}_S] ~  
\\
& \quad = ~ \Var' [(X) \mid \mathcal{F}_S] + 2 \Cov'[(X), (Y) \mid \mathcal{F}_S] + \Var'[(Y) \mid  \mathcal{F}_S]
& & \text{by Proposition~\ref{proposition:variance-sum}}
\\
& \quad = ~ \Var' [(X)\mid \mathcal{F}_S] & & \text{by Proposition~\ref{proposition:rpd-covariance}.}
\end{align*}
\end{proof}

\begin{corollary}
\label{corollary:variance-rate-preserved}
Suppose $(X)$ has variance-rate process $(A)$ and a path-integrable drift process $(B)$,
then $(A)$ is also the variance-rate process for the zero-drift process $(X-\int\! B)$  (see Corollary~\ref{corollary:kill-drift}).
\end{corollary}
\noindent
An alternative way of looking at Corollary~\ref{corollary:variance-rate-preserved} is that, when $(X)$ has 
path-integrable drift process $(B)$, variance rate can alternatively be defined as
a stopping derivative for 
\begin{equation}
\label{equation:nelson}
\textstyle
\Exp\left[\left. \left(X_{T-} - \left(\int \! B\right)_{T}  - X_S + \left(\int \! B\right)_{S}\,\right)^2 \,  \right|\,  \mathcal{F}_S  \right]\, \enspace ,
\end{equation}
which is the stopping functional defining relative-second-moment for the process $(X-\int\! B)$. 

We now have three different characterisations of variance rate as stopping derivatives: the original definition as a stopping derivative
for conditional variance~\eqref{equation:variance-rate}; the definition as relative-second-moment rate, i.e., as a stopping derivative for~\eqref{equation:rsm-rate}, which agrees with the previous for processes with drift; and the definition 
as a stopping derivative for~\eqref{equation:nelson}, which applies in the more restrictive case of processes with path-integrable drift processes. In all three cases, the stopping functionals utilised have a common format, returning the expected square distance from specified centres. In the case of conditional variance, the centre is the conditional expectation $\Exp[X_{T-} | \mathcal{F}_S]$. In the case of relative second moment, the centre is the starting value $X_S$. 
In the case of~\eqref{equation:nelson}, the centre is calculated according to  how the drift continuously varies between $S$ and $T$.
These are not the only possible choices of centre. 
Another natural choice is to use the drift $B_S := \Exp'[(X) \mid \mathcal{F}_S]$ at the start time $S$ to determine a projected centre $X_S + B_S (T-S)$ at time $T$, and 
hence consider a stopping derivative for
\begin{equation}
\label{equation:projected-centre}
\Exp[(X_{T-}  - X_S - B_S  (T-S))^2  \mid \mathcal{F}_S]\, \enspace .
\end{equation}
Indeed, in~\cite[\S5.1]{RWII}, an informal discussion of (what we are calling) variance rate motivates it as a form of right derivative for a similar formula. 
We use Proposition~\ref{proposition:rpd-covariance} to show that the stopping derivative for~\eqref{equation:projected-centre}  provides yet  another formula defining the variance rate. Curiously, this equivalence continues to hold if one uses any $\mathcal{F}_S$-measurable random variable whatsoever
for $B_S$ in~\eqref{equation:projected-centre}. It is  irrelevant whether or not
$B_S$ is the drift of $(X)$ at $S$.
\begin{proposition}
Suppose $(X)$ has drift at $S$ and  $B_S$ is any $\mathcal{F}_S$-measurable random variable.
Then the following are equivalent for
an $\mathcal{F}_S$-measurable random variable $C_S$.
\begin{enumerate}
\item $C_S$ is the variance rate for $(X)$ at $S$.
\item $C_S$ is a stopping derivative for~\eqref{equation:projected-centre}.
\end{enumerate}
\end{proposition}
\begin{proof}
We write $(B_S  (t-S))$ for the process $(B_S (t - S))_{t \geq S}$ defined as in~\eqref{equation:restrict-S}. This process is clearly right path-differentiable at $S$ with right path-derivative, hence stopping derivative, $B_S$. By Proposition~\ref{proposition:drift-sum}, the process $(X - B_S  (t-S))$ thus also has a stopping derivative
 at $S$, namely  $\Exp'[(X) \mid \mathcal{F}_S] -B_S$.

The stopping derivative for~\eqref{equation:projected-centre} is precisely the relative-second-moment rate for the process
$(X - B_S  (t-S))$, which by Proposition~\ref{proposition:rsm-rate} coincides with the variance rate for $(X - B_S  (t-S))$.
That is, the stopping derivative for~\eqref{equation:projected-centre} at $S$ is  equal to
$\Var'[(X - B_S  (t-S))  \mid \mathcal{F}_S]$, which is in turn equal to $\Var'[(X)  \mid \mathcal{F}_S]$ by
Proposition~\ref{proposition:variance-rate-preserved}.
\end{proof}


As a final application of covariance rate, we formulate the multi-dimensional It\^o formula. Accordingly, 
we consider processes valued in $\Real^n$, where for convenience we fix $n$. For such a process $(\VRV{X})$ we write 
$(\VRV{X}^1), \dots, (\VRV{X}^n)$ for the $n$ coordinate processes.  

\begin{definition}[Covariance-rate matrix] For a finite stopping time $S$, the $n \times n$ \emph{covariance-rate matrix}
$\COV'[(\VRV{X}) \mid \mathcal{F}_S]$ of $\mathcal{F}_S$-valued random variables is defined by:
 \[
\COV'[(\VRV{X}) \mid \mathcal{F}_S]_{i,j} ~ := ~ \Cov'[(\VRV{X}^i),(\VRV{X}^j) \mid \mathcal{F}_S] \enspace .
\]
\end{definition}

Proposition~\ref{proposition:non-negative} is the $n=1$ case of the more general result below. Although this is a basic result, we give the proof to highlight the measurability properties involved.

\begin{proposition}
\label{proposition:psd}
Whenever $\COV'[(\VRV{X}) \mid \mathcal{F}_S]$ is defined, it is $\mathcal{F}_S$-almost-surely positive semi-definite.
\end{proposition}

\begin{proof}
An $n \times n$ matrix $A$ is not positive semi-definite if and only if there exists an $n$-dimensional rational column vector $\mathbf{q}$ such that
$\mathbf{q}^\top \!A\, \mathbf{q} < 0$. Using this, if $\COV'[(\VRV{X}) \mid \mathcal{F}_S]$ is defined then 
the set 
\[W ~ := ~ \{\omega \in \Omega \mid \text{$\COV'[(\VRV{X}) \mid \mathcal{F}_S](\omega) $ is not positive semi-definite}\} \]
belongs to $\mathcal{F}_S$. Also, by systematically searching through rational $n$-vectors, there exists an $\mathcal{F}_S$-measurable
random $n$-vector $\VRV{Q}$ such that
\[
\VRV{Q}(\omega) ~:= ~ \begin{cases}
\text{$\mathbf{q}$ such that $\mathbf{q}^\top (\COV'[(\VRV{X}) \mid \mathcal{F}_S](\omega))\, \mathbf{q} < 0$} & \text{if $\omega \in W$} \\
\mathbf{0}& \text{otherwise.}
\end{cases}
\]
Suppose for contradiction that $W$ has positive probability. Define:
\[
\RVeps(\omega) ~ := ~ 
\begin{cases} -\frac{1}{2} (\VRV{Q}^\top \COV'[(\VRV{X}) \mid \mathcal{F}_S]\,  \VRV{Q}) (\omega) & \text{if $\omega \in W$} \\
1 & \text{otherwise.}
\end{cases}
\]
Clearly $\RVeps$ is $\mathcal{F}_S$-measurable and $\RVeps > 0$. We have:
\begin{align*}
& \VRV{Q}^\top \COV'[(\VRV{X}) \mid \mathcal{F}_S]\,  \VRV{Q}  ~ = ~ \sum_i \sum_j \VRV{Q}^i \Cov'[(\VRV{X}^i), (\VRV{X}^j) \mid \mathcal{F}_S] \,\VRV{Q}^j \enspace ,
\end{align*}
which, by the linearity of stopping derivatives, is the stopping derivative for
\begin{align*}
& 
\VRV{Q}^\top \COV[(\VRV{X}_{T-}) \mid \mathcal{F}_S]\,  \VRV{Q} \enspace .
\end{align*}
where $\COV[(\VRV{X}_{T-}) \mid \mathcal{F}_S]$ is the covariance matrix.  Accordingly there exists a finite stopping time $U > S$ such that, for every finite stopping time $T$ with $S \leq T \leq U$,
\begin{align}
\label{equation:sd-A}
 \VRV{Q}^\top \COV[(\VRV{X}_{T-}) \mid \mathcal{F}_S]\,  \VRV{Q} 
 ~ \leq ~
 (\VRV{Q}^\top \COV'[(\VRV{X}) \mid \mathcal{F}_S]\,  \VRV{Q} + \RVeps) \, 
 \Exp[T-S  \mid \mathcal{F}_S] \enspace .
 \end{align}
 In the case $T := U$, we have $ \VRV{Q}^\top \COV[(\VRV{X}_{T-}) \mid \mathcal{F}_S]\,  \mathbf{Q}  \geq 0$,
 since covariance matrices are positive semi-definite.
 However, also $((\mathbf{Q}^\top \COV'[(\VRV{X}) \mid \mathcal{F}_S]\,  \mathbf{Q} + \RVeps) \, 
 \Exp[T-S  \mid \mathcal{F}_S])(\omega) < 0$, for all $\omega \in W$, 
 by the choice  of $\RVeps$. This contradicts~\eqref{equation:sd-A} above.
\end{proof}

\begin{theorem}[The mult-dimensional It\^o formula]
\label{theorem:multi-dimensional-Ito}
Suppose that  $(\VRV{X})$ is an $n$-dimensional cadlag process with drift and covariance-rate matrix at a finite stopping time $S$.
Suppose that $\underline{f} \colon \Omega \times \Real^n \to \Real$ is an  $\mathcal{F}_S$-jointly-measurable random function such that
$\underline{f}_\omega(\VRV{X}_t(\omega))$ is defined for all $t \geq S(\omega)$.
Suppose further that  $\underline{f}$ is twice
$\mathcal{F}_S$-measurably continuously differentiable at $\VRV{X}_S$.
Then the process $(\underline{f}(\VRV{X}))$ has both drift 
and variance rate at $S$, with these given by:
\begin{align}
\nonumber
\Exp'[ (\underline{f}(\VRV{X})) \mid \mathcal{F}_S] ~ = & ~  \sum_{i=1}^{n} \partial_i \underline{f}\,(\VRV{X}_S)\,\Exp'[(\VRV{X}^i) \mid \mathcal{F}_S]
\\ 
\label{equation:multi-ito-exp}
& \qquad  +
 \frac{1}{2}\sum_{i=1}^{n} \sum_{j=1}^{n} \partial_i\partial_j\underline{f}\,(\VRV{X}_S)\,\Cov'[(\VRV{X}^i),  (\VRV{X}^j)\mid \mathcal{F}_S] \enspace , \\
 \label{equation:multi-ito-var}
\Var'[ (\underline{f}(\VRV{X})) \mid \mathcal{F}_S] ~ = & ~ \sum_{i=1}^{n} \sum_{j=1}^{n}\,
(\partial_i\underline{f}\,(\VRV{X}_S))(\partial_j\underline{f}\,(\VRV{X}_S)) \,\Cov'[(\VRV{X}^i),  (\VRV{X}^j) \mid \mathcal{F}_S]\enspace .
\end{align}
\end{theorem}

Following a time-honoured tradition, we omit the proof of Theorem~\ref{theorem:multi-dimensional-Ito}, which closely follows the proof of Theorem~\ref{theorem:ito-formula}, but with notational overhead due to the multi-dimensional setting. One further remark is perhaps helpful. In many places in the proof of  Theorem~\ref{theorem:multi-dimensional-Ito},
the nonnegative term 
$(X_{T\!-}(\omega) \!- \!X_{S}(\omega))^2$  in the proof of Theorem~\ref{theorem:ito-formula} is replaced with 
\begin{equation}
\label{eq:non-polarised}
(\VRV{X}^i_{T\!-}(\omega) \!- \!\VRV{X}^i_{S}(\omega))(\VRV{X}^j_{T\!-}(\omega) \!- \!\VRV{X}^j_{S}(\omega)) \enspace ,
\end{equation}
which may assume arbitrary real values. This complicates the process of finding 
 lower and upper bounds in the analogue of~\eqref{equation:AB-y}. 
Such bounds are  obtained via the polarisation formula
\[
\frac{1}{4}
((\VRV{X}^i_{T\!-}(\omega) + \VRV{X}^j_{T\!-}(\omega)) - (\VRV{X}^i_{S}(\omega) + \VRV{X}^j_{S}(\omega))^2
-
 \frac{1}{4}((\VRV{X}^i_{T\!-}(\omega) - \VRV{X}^j_{T\!-}(\omega)) - (\VRV{X}^i_{S}(\omega) - \VRV{X}^j_{S}(\omega))^2 
\enspace ,
\]
which expresses~\eqref{eq:non-polarised} as the difference of two nonnegative expressions.


\section{Local martingales}
\label{section:fundamental}

In this section, we prove that continuous zero-drift processes coincide with random translations of continuous local martingales.
%
This result will provide us with a bridge between the stopping-derivative-based methods of the present article and the 
It\^o-integral-based methods of stochastic calculus in its usual guise, which are heavily based on local martingales.
\begin{lemma}
\label{lemma:driftless}
If $(X)$ is a zero-drift process and $T$ is a stopping time, then the stopped process
$(X^T)$ also has zero drift.
\end{lemma}

\noindent
We omit the easy proof.

\begin{theorem} \leavevmode
\label{theorem:driftless-martingale}
The following are equivalent for a continuous process $(X)$.
\begin{enumerate}
\item \label{dm:i} $(X)$  has zero drift.
\item \label{dm:ii} $(X - X_0)$ is a local martingale.
\item \label{dm:iii} $(X)$ can be written as $(Y + Z_0)$, where $(Y)$ is a continuous local martingale and $Z_0$ is an $\mathcal{F}_0$-measurable random variable.
\end{enumerate}
\end{theorem}
\begin{proof}
%
%
The implication \ref{dm:ii}~$\implies$~\ref{dm:iii} is trivial, and the implication~\ref{dm:iii}~$\implies$~\ref{dm:i} follows from Theorem~\ref{theorem:martingale-driftless} and elementary properties of drift.

To prove \ref{dm:i}~$\implies$~\ref{dm:ii}, 
we claim that every  continuous bounded zero-drift process $(X)$ with $X_0 = 0$ is a martingale. 
The desired implication  follows, because if $(Y)$ is a continuous zero-drift process, then let $H_n$ be the debut time
for $(Y- Y_0)$ landing in the set $\{x \in \Real \mid |x| \geq n\}$. For every $n \in \Nat$, the stopped process $((Y-Y_0)^{H_n})$ is  continuous, bounded and, by Lemma~\ref{lemma:driftless}, has zero drift. Thus, by the claim, every $((Y-Y_0)^{H_n})$ is a martingale. 
Since $(Y-Y_0)$ is continuous, $(H_n) \uparrow \infty$. Thus $(H_n)$ provides a sequence of stopping times showing that $(Y-Y_0)$ is a local martingale.

It remains to prove the claim. Accordingly, suppose that $(X)$ is a continuous bounded zero-drift process with $X_0 = 0$.
For any $s < t$ and $\varepsilon > 0$, we shall prove that 
\begin{equation}
\label{epsilon-martingale}
|\Exp[X_t | \mathcal{F}_s] - X_s| \leq \varepsilon(t-s) \enspace .
\end{equation}
It follows that
$\Exp[X_t | \mathcal{F}_s]  = X_s$ for all $s < t$; i.e., $(X)$ is indeed a martingale. 

For the proof of~\eqref{epsilon-martingale}, we  
construct a transfinite sequence $(S_\alpha)_{\alpha < \omega_1}$ (where $\omega_1$ is the smallest uncountable ordinal) of stopping times as follows.
Define $S_0 := s$. Given $S_\alpha$, define $S_{\alpha + 1} := U_{\alpha  + 1} \wedge t$, where $U_{\alpha+1}$ is a chosen stopping time, as given by the zero-drift property, such
that $S_\alpha < U_{\alpha+1}$ and,  for all $T$ with $S_\alpha \leq T \leq U_{\alpha+1}$, 
\begin{equation}
\label{def-by-tr}
|\Exp[X_T - X_{S_\alpha} \mid \mathcal{F}_{S_\alpha}]| \leq  \varepsilon \, \Exp[T-S_\alpha  \mid \mathcal{F}_{S_\alpha}] \enspace .
\end{equation}
 For a limit ordinal $\lambda < \omega_1$, define
$S_\lambda = \sup_{\alpha < \lambda} S_\alpha$. Then, for $\alpha \leq \beta < \omega_1$, we have $s \leq S_\alpha \leq S_\beta \leq t$. 

We prove by transfinite induction that, for all $\alpha < \omega_1$, 
\begin{equation}
\label{prove-by-ti}
|\Exp[X_{S_\alpha} | \mathcal{F}_s] - X_s| \leq \varepsilon  \, \Exp[S_\alpha -s \mid \mathcal{F}_s ]\enspace .
\end{equation}
The case $\alpha = 0$ is immediate because $S_0 = s$.
For successor ordinals, given~\eqref{prove-by-ti}  as induction hypothesis and using the fact that $S_\alpha \leq S_{\alpha+1} \leq U_{\alpha+1}$, we have:
\begin{align*}
 & |\Exp[X_{S_{\alpha +1}} | \mathcal{F}_s] - X_s|   \\
 & \quad = 
|\Exp[(\Exp[X_{S_{\alpha +1}} | \mathcal{F}_{S_\alpha} ] - X_{S_\alpha} )| \mathcal{F}_s] +
\Exp[X_{S_{\alpha}} | \mathcal{F}_s] - X_s |
\\
& \quad \leq
\Exp[(|\Exp[X_{S_{\alpha +1}} | \mathcal{F}_{S_\alpha} ] - X_{S_\alpha} |)| \mathcal{F}_s] +
|\Exp[X_{S_{\alpha}} | \mathcal{F}_s] - X_s |
\\
& \quad \leq \Exp[\varepsilon \, \Exp[S_{\alpha+1}-S_\alpha  \mid \mathcal{F}_{S_\alpha}] \mid  \mathcal{F}_s] +
\varepsilon\,\Exp[S_\alpha -s \, | \mathcal{F}_s]  & & \text{by \eqref{def-by-tr} and \eqref{prove-by-ti}}\\
&  \quad = \varepsilon (\Exp[S_{\alpha+1}-S_\alpha \, | \mathcal{F}_s] +
\Exp[S_\alpha -s \, | \mathcal{F}_s])
\\
&  \quad = \varepsilon\, \Exp[S_{\alpha+1}-s \, | \mathcal{F}_s] \enspace .
\end{align*} 
For a limit ordinal $0 < \lambda < \omega_1$, let $(\alpha_n)$ be an increasing sequence of smaller ordinals such that
$\lambda = \sup_n \alpha_n$. Then $(S_{\alpha_n}) \uparrow S_\lambda$. So, since $X$ has continuous paths,
it follows that $(X_{S_{\alpha_n}}) \to X_{S_\lambda}$ converges almost surely. Because the process $(X)$ is bounded,
\begin{equation}
\label{X-convergence}
(\Exp[X_{S_{\alpha_n}} | \mathcal{F}_s]) \to \Exp[X_{S_\lambda} |\mathcal{F}_s]
\end{equation}
converges almost surely, by 
the dominated convergence theorem. By induction hypothesis, for every $\alpha_n$, we have
\[
|\Exp[X_{S_{\alpha_n}} | \mathcal{F}_s] - X_s| \leq \varepsilon  \, \Exp[S_{\alpha_n} - s \mid  \mathcal{F}_s ]
\leq \varepsilon  \, \Exp[S_\lambda - s \mid  \mathcal{F}_s ] \enspace .\]
It thus follows from~\eqref{X-convergence} that
$|\Exp[X_{S_{\lambda}} | \mathcal{F}_s] - X_s| 
\leq \varepsilon  \, \Exp[S_\lambda - s \mid  \mathcal{F}_s ]$, as required. This completes the proof of~\eqref{prove-by-ti}.

We are now in a position to prove~\eqref{epsilon-martingale}. Note that
if $S_\alpha < t$ holds on a positive measure subset $W \subseteq \Omega$, then $S_\alpha < S_{\alpha+1}$ on $W$, hence
$\Exp[S_{\alpha + 1}] -\Exp[S_\alpha] = \Exp[S_{\alpha + 1} - S_{\alpha}] > 0$.

Let $O$ be the set of ordinals $\alpha < \omega_1$ such that $S_\alpha < t$ on some positive measure subset, equivalently such that
$\Exp[S_\alpha] < t$. Clearly $O$ is a down-closed set of countable ordinals. By the inequality at the end of the last paragraph, $\alpha < \beta \in O$ implies
$\Exp[S_\alpha] < \Exp[S_\beta]$; i.e., the map $\alpha \mapsto \Exp[S_\alpha] : O \to [s,t]$ is strictly increasing. Since there is no strictly increasing map into $\Real$ whose domain is the set of all countable ordinals, there exists a smallest countable ordinal $\gamma$ such that $\gamma \notin O$, that is, such that $\Exp[S_\gamma] = t$. We thus have 
$S_\gamma = t$ almost surely. 

Using~\eqref{prove-by-ti} for the inequality, it follows that 
\[|\Exp[X_t | \mathcal{F}_s] - X_s| = |\Exp[X_{S_\gamma} | \mathcal{F}_s] - X_s| \leq \varepsilon \,\Exp[S_\gamma -s] = \varepsilon (t-s) \enspace . \]
This establishes~\eqref{epsilon-martingale}, completing the proof of the theorem.
\end{proof}

\noindent
Property~\ref{dm:iii} of Theorem~\ref{theorem:driftless-martingale} defines what we mean by \emph{random translations of local martingales}. Such random translations take one outside the class of local martingales precisely in cases in which the displacement random variable  is not integrable. 
\begin{corollary}
\label{corollary:driftless-martingale}
The following are equivalent for a continuous process $(X)$.
\begin{enumerate}
\item \label{cdm:i} $(X)$  has zero drift and $X_0$ is integrable.
\item \label{cdm:ii} $(X)$ is a local martingale.
\end{enumerate}
\end{corollary}
\begin{proof}
If property~\ref{cdm:i} holds then, by Theorem~\ref{theorem:driftless-martingale}, $(X - X_0)$ is a local martingale, hence so is $(X)$, since local martingales are preserved under translation by
integrable $\mathcal{F}_0$-measurable random variables.

Conversely, if $(X)$ is a local martingale then it has zero drift by Theorem~\ref{theorem:martingale-driftless}, and the integrability of $X_0$ holds because there exist finite stopping times $T>0$ for which the stopped process
$(X^T)$ is a martingale, which means $X^T_0$ is integrable, and we have $X_0 = X^T_0$.
\end{proof}

\begin{remark}
\label{remark:non-martingale-example}
The relationship between zero drift and the local martingale property does not extend to cadlag processes. A simple example of a (bounded) cadlag zero-drift process that is not a local martingale is the deterministic process $(x_t)_t$ defined by:
\[ x_t ~ = ~ 
\begin{cases} 
0 & \text{if $t = 0$,} \\
n^{-2}& \text{if $n \in  \mathbb{Z}$,  $n > 0$ and  $n^{-1}  \leq t < (n-1)^{-1}$,}
\end{cases}
\]
where $0^{-1}$ is interpreted as $\infty$. 
I do not know any example of a cadlag local martingale that does not have zero drift. 
\end{remark}

\begin{corollary}
\label{corollary:unique}
The zero process $(0)$ is the only continuous process $(X)$ with drift process $(0)$ and variance-rate process $(0)$ satisfying $X_0 = 0$. 
\end{corollary}
\noindent
Note that the example from Remark~\ref{remark:non-martingale-example} shows that the continuity assumption cannot be weakened to cadlag.
\begin{proof}
Suppose that $(X)$ is a continuous process with drift process $(0)$ and variance-rate process $(0)$ and such that $X_0 =0$. 
By~\eqref{equation:var-square}, we have   $ (\Exp'[(X^2)  \mid \mathcal{F}] ) = (\Var'[(X)   \mid  \mathcal{F}]) = (0)$; i.e.,
the process $(X^2)$ has zero drift. Since $(X^2)$ is continuous, it is, by Corollary~\ref{corollary:driftless-martingale}, a local martingale. It is, furthermore,  a non-negative local martingale equal to $0$ at time $0$. 
An easy argument shows that the zero process $(0)$ is the only local martingale with these properties. So $(X^2) = (0)$, whence $(X) = (0)$.
\end{proof}

We say that a process $(X)$ has \emph{zero variance rate} if it has drift and its variance rate process is $(0)$. By Proposition~\ref{proposition:rpd-covariance}, every right-path-differentiable process has zero variance rate. The next result strengthens Corollary~\ref{corollary:kill-drift} with a uniqueness 
property. 
\begin{proposition}
\label{proposition:unique-int-B}
Suppose $(X)$ has a path-integrable drift process $(B)$, then there exists a  unique continuous zero-variance-rate process $H$ with $H_0=0$ such that
$(X-H)$ is a zero-drift process, namely  $(H) := (\int\! B)$.
\end{proposition}
\begin{proof}
Firstly, the process $(\int\! B)$ is indeed continuous, zero at time $0$, and has zero variance rate by Propositions~\ref{proposition:integrable-gives-drift} and 
\ref{proposition:rpd-covariance}. Moreover, $(X-\int\! B)$ has zero drift by Corollary~\ref{corollary:kill-drift}.

For uniqueness, suppose $(H)$ is any continuous zero-variance-rate process such that $H_0 = 0$ and $(X-H)$ has zero drift.
By the sum rule for drift (Proposition~\ref{proposition:drift-sum}), the process $({\int\! B} - H)$ has zero drift. It is also clearly continuous and zero at time $0$.
By Proposition~\ref{proposition:variance-rate-preserved},
\[ \textstyle
\Var'[({\int\! B} - H) \mid \mathcal{F}_S] ~ =  ~ \Var'[(H) \mid \mathcal{F}_S]  
~ = ~ 0 
\enspace .\]
Hence it follows from Corollary~\ref{corollary:unique} that
$({\int\! B} - H) = (0)$, i.e., $({\int\! B}) = ( H)$.
\end{proof}

\begin{proposition}
\label{proposition:unique-decomposition}
Suppose $(X)$ has both drift process and variance-rate process with the former path integrable. Then $(X)$ has a unique decomposition $(X) = (Y) + (Z)$ where $(Y)$ has zero variance rate with $Y_0 = 0$ and the $(Z)$ has zero drift. Moreover,
\[
\Exp'[(Y) \mid \mathcal{F}] ~ = ~ \Exp'[(X) \mid \mathcal{F}] \quad \text{and} \quad 
\Var'[(Z) \mid \mathcal{F}] ~ = ~ \Var'[(X) \mid \mathcal{F}]  \enspace .
\]
\end{proposition}
\begin{proof}
By Proposition~\ref{proposition:unique-int-B}, defining $(Y) := (\int\! B)$, we obtain that
$(Y)$ is the unique zero-variance-rate with process with $Y_0 = 0$ such that  and $(Z) := (X - Y)$ has zero drift.
Clearly $\Exp'[(Y) \mid \mathcal{F}] = \Exp'[(X) \mid \mathcal{F}]$. Moreover, 
$\Var'[(Z) \mid \mathcal{F}] = \Var'[(X) \mid \mathcal{F}]$ by Corollary~\ref{corollary:variance-rate-preserved}.
\end{proof}

\begin{proposition}
\label{proposition:unique-int-A}
Suppose $(X)$ is a continuous zero-drift process. Then the following are equivalent for a path-integrable 
process $(A)$.
\begin{enumerate}
\item $(A) = \Var'[(X) \mid \mathcal{F}]$.
\item \label{item:qv} $(X^2-\int\! A )$ has zero drift.
\end{enumerate}
Moreover, when $(A) = \Var'[(X) \mid \mathcal{F}]$ 
then $([X]) := (\int\! A)$ is the unique continuous, zero-variance-rate, initially zero process for which
$(X^2-[X])$ has zero drift. Furthermore, $([X])$ is increasing.
\end{proposition}
\begin{proof}
%
Suppose that $(X^2-\int\! A )$ has zero drift. Then
\[\textstyle
(\Var'[(X)   \mid  \mathcal{F}]) ~ = ~ (\Exp'[(X^2)  \mid \mathcal{F}] ) ~ = ~
(\Exp'[(\int \! A)  \mid \mathcal{F}] )~ = ~ (A) \enspace ,\]
where  the first equality is~\eqref{equation:var-square} and the last holds by Proposition~\ref{proposition:path-differentiable}.

Conversely, suppose $(A) = (\Var'[(X)   \mid  \mathcal{F}])$. Since $(\Var'[(X)   \mid  \mathcal{F}]) = (\Exp'[(X^2)  \mid \mathcal{F}] )$ by~\eqref{equation:var-square}, 
it is immediate
from Proposition~\ref{proposition:unique-int-B} that
$([X]) := (\int\! A)$ is the unique continuous, zero-variance-rate, initially zero process for which
$(X^2-[X])$ has zero drift. Moreover, $([X])$ is increasing because it is defined as $(\int\! A)$, where $(A)$ is non-negative by Proposition~\ref{proposition:non-negative}.
\end{proof}
\noindent
By  the equivalence between zero-drift processes and local martingales, 
in the case that ${X_0}^2$ is integrable, 
property~\ref{item:qv} of the above proposition asserts that $(\int\!A)$ satisfies the characterising property of 
the \emph{quadratic variation process} for $(X)$. So we have a  simple construction of the quadratic variation process that, although not applicable to all continuous local martingales, does apply to those with  path-integrable variance rate. Since this holds for processes defined as stochastic integrals (see
Theorem~\ref{theorem:stochastic-integral} below for a precise formulation), this is still a very broad class.

As a further consequence of the equivalence between continuous zero-drift processes 
and continuous local martingales, we obtain a version of 
Levy's characterisation of Brownian motion formulated in terms of drift and 
variance rate.
\begin{theorem}
\label{theorem:levy}
The following are equivalent for a continuous process $(X)$ with $X_0 = 0$.
\begin{enumerate}
\item \label{levy:a} $X$ has drift process $(0)$ and variance-rate process $(1)$.
\item \label{levy:b} $X$ is a standard Brownian motion.
\end{enumerate}
\end{theorem}
\begin{proof} The implication~\eqref{levy:b} $\Rightarrow$~\eqref{levy:a} is Proposition~\ref{proposition:brownian}.

The quick proof of~\eqref{levy:a} $\Rightarrow$~\eqref{levy:b} 
first uses Theorem~\ref{theorem:driftless-martingale}
to translate 
having drift process $(0)$ to $(X)$ being a local martingale, and 
having variance rate $(1)$ to $({X_t}^2 - t)_t$ being a local martingale. It then follows that $(X)$ is a standard Brownian motion by Levy's theorem in its usual modern local-martingale formulation.
\end{proof}
\begin{remark}
The above short argument depends on the local-martingale formulation of Levy's theorem, which is often proved using the machinery of It\^o-integral-based stochastic calculus. Actually this machinery can be bypassed entirely, by taking the following alternative route. 
Instead of applying Theorem~\ref{theorem:driftless-martingale}
at the beginning of the argument, one can follow the  elegant Kunita-Watanabe proof of Levy's theorem, which can be recast entirely in terms of our
stopping-derivative-based version of the Itô formula (Theorem~\ref{theorem:ito-formula}),
making direct use of the  drift process $(0)$ and variance rate $(1)$.
In this way, one obtains a zero-drift process involving characteristic functions, which one needs to show is a martingale.
Only at this point is Theorem~\ref{theorem:driftless-martingale} applied, yielding that this process is a local martingale, after which the usual argument that it is actually a martingale applies. 
\end{remark}

\section{Stochastic integration and uniqueness}
\label{section:integration-uniqueness}

Up to this point, our results have been obtained without use of  stochastic integration.
As in ordinary calculus, one of the  basic motivations for integration is that it  provides a mechanism
for constructing processes with specified derivatives. In our case the relevant derivatives are the  drift and variance rate, defined as stopping derivatives, and the relevant form of integration is the stochastic integral of It\^o.
Accordingly, suppose $(B)$ and $(A)$ are path-integrable $(\mathcal{F})$-adapted processes with $(A)$ non-negative. Assuming the filtration $(\mathcal{F})$ is rich enough to support a Brownian motion, 
Theorem~\ref{theorem:stochastic-integral} below shows that it is a simple application of stochastic integration to 
construct a continuous process $(X)$ such that
\begin{equation}
\label{eq:specification}
\Exp'[(X) \mid \mathcal{F}] ~ = ~ (B)\quad \text{and} \quad 
\Var'[(X) \mid \mathcal{F}] ~ = ~ (A)  \enspace .
\end{equation}

As a first step, we characterise \eqref{eq:specification} above as equivalent to the process $(X)$ solving a  {local martingale problem}.
In this and henceforth, 
we use the equivalence between zero drift and the local martingale property
(Corollary~\ref{corollary:driftless-martingale})
without explicit mention.

\begin{proposition}
\label{proposition:local-martingale-problem}
Suppose $(A)$ and $(B)$ are path-integrable, then the following are equivalent for a continuous process $(X)$
with integrable $X_0$.
\begin{enumerate}
\item \label{item:drift-plus-variance-rate} 
$(X)$ has drift process $(B)$ and variance rate process $(A)$ (i.e., \eqref{eq:specification} above).
\item \label{item:local-martingale-problem}
$\left(X - {\int \! B}\right)$ is a local martingale with quadratic variation process $\left({\int \! A}\right)$.
\end{enumerate}
\end{proposition}
\begin{proof}
Suppose $\left(X - {\int \! B}\right)$ is a local martingale with quadratic variation process $\left({\int \! A}\right)$.
Then $\Exp'[X - {\int \! B} \mid \mathcal{F}] = (0)$, so $\Exp'[(X) \mid \mathcal{F}] = (B)$, by 
Proposition~\ref{proposition:integrable-gives-drift}. The defining property of the quadratic variation process gives
that $\left(\left(X - {\int \! B}\right)^2 - {\int \! A} \right)$ is a local martingale.
So $\Var'\left[\left.\left(X - {\int \! B} \right) \right| \mathcal{F} \right] = (A)$, by Proposition~\ref{proposition:unique-int-A}.
Hence $\Var'\left[\left.\left(X\right) \right| \mathcal{F} \right] = (A)$, by Corollary~\ref{corollary:variance-rate-preserved}.

Conversely, suppose $(X)$ has drift process $(B)$ and variance rate process $(A)$. 
Then $\left(X - {\int \! B}\right)$ is a local martingale by Corollary~\ref{corollary:kill-drift}.
Also 
\[\textstyle \Exp'\left[\left.\left(\left(X - {\int \! B}\right)^2\right)\right| \mathcal{F} \right] ~ = ~
\Var'\left[\left.\left(X - {\int \! B}\right)\right| \mathcal{F} \right] ~ = ~
\Var'\left[\left.\left(X\right)\right| \mathcal{F} \right]  ~ =  ~ (A) \enspace ,
\]
by \eqref{equation:var-square} and Corollary~\ref{corollary:variance-rate-preserved}. 
So $\left(\left(X - {\int \! B}\right)^2 - \int \! A\right)$ is a local martingale by Corollary~\ref{corollary:kill-drift}.
That is, $\left(\int \! A\right)$ is the quadratic variation process for $\left(X - {\int \! B}\right)$.
\end{proof}

\begin{theorem}[Fundamental theorem of calculus for stopping derivatives]
\label{theorem:stochastic-integral}
Suppose there exists some Brownian motion $(W)$ relative to $(\mathcal{F})$, 
suppose $C_0$ is $\mathcal{F}_0$-measurable, and 
$(\sigma^2)$ and $(B)$ are 
path integrable. Then 
\[
X_t  ~ := ~ C_0 ~ + ~ \int_0^t B_s \, ds  ~ + ~ \int_0^t \sigma_s \, dW_s
\]
defines a continuous process $(X)$ 
with drift process $(B)$ and variance-rate process $(\sigma^2)$ satisfying $X_0 = C_0$.
\end{theorem}
\begin{proof}
It is standard that the stochastic integral $Y_t := \int_0^t {\sigma_s} \, dW_s$ defines a continuous local martingale $(Y)$
with quadratic variation process ${\int\! \sigma^2}$. Since $(Y) = (X - C_0 - \int\! B)$, it is 
immediate from Proposition~\ref{proposition:local-martingale-problem} that 
$(X-C_0)$ 
has drift process $(B)$ and variance-rate process $(\sigma^2)$. Hence $(X)$ does too. 
%
%
\end{proof}
\noindent
We remark that Theorem~\ref{theorem:stochastic-integral} includes, as a special case, the following
global version 
of the right-differentiable formulation of the deterministic fundamental theorem of calculus:
if a right-continuous function $f$ on $[0,\infty)$ satisfies $\int_0^t |f(s)| {\mathit{d}s} < \infty$, for all $t$,  then the continuous function $t \mapsto \int_0^t f(s) {\mathit{d}s}$ is right differentiable with right derivative $f$. 
In the case of Theorem~\ref{theorem:stochastic-integral}, the right-continuity assumption, which is clearly needed, is built into our definition of path-integrability (Definition~\ref{definition:path-integrable}).

Theorem~\ref{theorem:stochastic-integral} recovers $(\sigma^2)$ as the variance rate, but does not  directly recover $(\sigma)$ itself. This is by necessity, as the $\sigma_s$ values are  not quantities that are intrinsic to the process $(X)$. They rather 
expresses a relationship between $(X)$ and the background Brownian motion $(W)$. We remark that~\cite{Allouba} presents a stochastic fundamental theorems of calculus that  does recover $\sigma_s$ itself, by making use  of an  exotic form of time derivative at $s$
for the quadratic covariation of $(X)$ and $(W)$.  
A more systematic result in a similar direction appears as \cite[Theorem~5.8]{CF},
where the process $(\sigma)$
is recovered as the 
 \emph{vertical derivative} $\nabla_{W}X$ of $(X)$ with respect to $(W)$, defined 
using  the  \emph{functional It\^{o} calculus}, whose development is the main contribution of
\emph{op.\ cit.}

Theorem~\ref{theorem:stochastic-integral} gives an existence result: given $C_0$, $(B)$ and nonnegative $(A) = (\sigma^2)$, one can find a process $(X)$ that has these as initial distribution, drift process and variance-rate process respectively. Since one could replace $\sigma$ in Theorem~\ref{theorem:stochastic-integral}
with $-{\sigma}$, it is not surprising that $C_0$, $(B)$ and $(\sigma^2)$  do not characterise $(X)$, not even up to equality in distribution, as the example below illustrates. Recall that two continuous processes $(X)$ and $(Y)$ are said to be \emph{equal in distribution} if $(X)$ and $(Y)$ induce the same probability measure on the standard Borel space $C([0,\infty))$ of continuous functions from $[0,\infty)$ to $\mathbb{R}$.

\begin{example}
\label{example:non-unique}
Let $(W)$ be a standard Brownian motion with respect to a filtration $(\mathcal{F})$. Consider the
stopping time 
\[ T(\omega) ~ := ~ \inf\{t \mid W_{t}(\omega) \geq 1\} \enspace .
\]
Define $(B) := (0)$ and $A_t := \One_{t < T}$. Then the stopped process $(W^T)$ and its negation
$(-W^T)$ both have initial value $0$, drift process $(B)$ and variance-rate process $(A)$. 
However, they are not equal in distribution since $\lim_{t \to \infty} P(W^T_t = 1) = 1$, whereas
$\lim_{t \to \infty} P(-W^T_t = 1) = 0$.
\end{example}

We end this section with a uniqueness result, 
chosen to showcase our stopping-derivative-based calculus:  the result has an  easy  proof using tools we already have at our disposal.
We give
further conditions under which a continuous local martingale $(X)$ with $X_0 = 0$ is determined up to equality-in-distribution by its variance-rate process $(A)$.
The main idea is to require $A_t$ 
to be functionally dependent on the trajectory of $(X)$ over time $[0,t]$.
Under certain conditions, the process $(X)$ is determined up to equality-in-distribution by the function underlying this dependency,
which we call a  \emph{$C_0([0,\infty))$-adaptation}, where 
\[C_0([0,\infty)) := \{f \in C([0,\infty)) \mid f(0) = 0\}\, . \]

\begin{definition}[$C_0([0,\infty))$-adaptation]
A \emph{$C_0([0,\infty))$-adaptation} is a function
\[
a \colon C_0([0,\infty)) \times [0,\infty) \to \Real 
\]
satisfying:
\begin{enumerate}
\item for all $t \geq 0$, the function $f \mapsto a(f,t) \colon C_0([0,\infty)) \to \Real$ is measurable, and
\item for all $f,g \in C_0([0,\infty))$ and $t \geq 0$, if $f\!\restriction_{[0,t]} = g \! \restriction_{[0,t]}$.
then $a(f,t) = a(g,t)$, and
\end{enumerate}
We say that $a$ is \emph{right-continuous} if further
\begin{enumerate}
\setcounter{enumi}{2}
\item for all $f \in C_0([0,\infty))$, the function $t \mapsto a(f,t) \colon [0,\infty) \to \Real$ is right-continuous.
\end{enumerate}
\end{definition}

\noindent
Given a continuous process $(X)$, any $C_0([0,\infty))$-adaptation $a$ induces
a process $a[(X)] := (a((s \mapsto X_s), t))_{t \geq 0}$, which is right-continuous whenever the adaptation $a$ is.
For instance, in Example~\ref{example:non-unique}, it holds that $A = a[(W)]$, where $a$ is the
right-continuous $C_0([0,\infty))$-adaptation
\[
a(f,t) ~ := ~ \begin{cases}
 1 & \text{if $\max \{f(s) \mid 0 \leq s \leq t\} < 1$,} \\
 0 & \text{otherwise.}
 \end{cases}
\]
It similarly holds that $A = a'[(-W)]$, for an evident  right-continuous $C_0([0,\infty))$-adaptation
$a'$ that is different from $a$.

Theorem~\ref{theorem:unique} below shows that, when 
a $C_0([0,\infty))$-adaptation exists that induces the variance-rate process of 
a continuous local martingale $(X)$ with $X_0 = 0$, then, 
as long as certain further conditions are satisfied,  the process $(X)$
is determined up to equality-in-distribution by the adaptation function.
The further conditions are that $a$ is positive and has \emph{unbounded path integrals}, where the latter means that, 
for every $f \in C_0([0,\infty))$ and  $t < \infty$, 
\[
\int_0^t a(f,s)\, ds < \infty ~~~ \text{and} ~~~
\int_0^\infty a(f,s)\, ds  = \infty \enspace .
\]

\begin{lemma}
Suppose  $a$  is a positive right-continuous $C_0([0,\infty))$-adaptation with unbounded path integrals.
Then the function $\Phi_a \colon C_0([0,\infty)) \to C_0([0,\infty))$ defined by
\begin{equation}
\label{equation:Phi}
\Phi_a(f) ~ := ~ s \mapsto f\!\left(\inf\left\{ t \, \left| \, \int_0^t a(f,u)\, du \geq s \right . \right\} \right)
\end{equation}
is a measurable bijection, hence has measurable inverse.
\end{lemma}

\noindent
We omit the proof, which is an exercise in real analysis.

\begin{theorem}
\label{theorem:unique}
Suppose $(X)$ is a continuous local martingale relative to $(\mathcal{F})$ with $X_0 = 0$.
Suppose further that there exists a  right-continuous $C_0([0,\infty))$-adaptation $a$ such that
$(X)$ has 
variance-rate process $a[(X)]$, 
where $a$ is positive and has unbounded path integrals.
Then the probability law $P_{(X)}$ of the process $(X)$ is given by the pushforward $(\Phi_a^{-1})_*(\mu_W)$  of the Wiener measure
$\mu_W$ on $C_0([0,\infty))$ along 
$\Phi_a^{-1}$, where $\Phi_a$ is as defined in \eqref{equation:Phi}.
\end{theorem}
\begin{proof}
We first obtain a Brownian motion $(W)$ as a time
change of $(X)$, after which we recover $(X)$ from $(W)$ using $\Phi_a^{-1}$.

Define a time-change $(R_s)_{s \geq 0}$ relative to $(\mathcal{F})$ by
\[
R_s(\omega)~ := ~ \inf \left\{ t \, \left| \,  \int_0^t a[(X)]_u(\omega)  \,du \geq  s   \right.\right\} \enspace ,
\]
calculating the time $t$ needed for the integral of the trajectory of $a[(X)]$ to reach $s$. Define the path-integral process
$(R') := (\int\!a[(X)])$, which is obviously adapted to $(\mathcal{F})$, and is also a time-change relative to $(\mathcal{F}_R)$.
Since $a$ is positive with $\int_0^\infty a(f,s)\,ds = \infty$ for all $f \in C([0,\infty))$, for every $\omega$, we have
$t \mapsto R'_t(\omega)$ is an increasing right-differentiable bijection from $[0,\infty)$ to itself with inverse
$s \mapsto R_s(\omega)$. Clearly $(D') := (a[(X)])$ is the
right-path-derivative process for $(R')$. Since every trajectory $t \mapsto R'_t(\omega)$ of $(R')$ is an increasing right-differentiable bijection with
positive right derivative, its inverse $s \mapsto R_s(\omega)$ is also right-differentiable with right derivative
\[D_s ~ = ~ \frac{1}{a[(X)]_{R_s}} \enspace .\]
This defines an $(\mathcal{F}_R)$-adapted  right-continuous process $(D)$ which is a right-derivative-process for $(R)$.

Define the time-changed process $(W) := (X_R)$. By Proposition~\ref{proposition:time-change}, for any finite 
$(\mathcal{F}_R)$ stopping time $S$, we have
\[
\Exp'[(W_S) \mid \mathcal{F}_{R_S}] ~ = ~ \Exp'[(X) \mid \mathcal{F}_{R_S}] \cdot D_S ~ = ~ 0 \enspace ,
\]
because $(X)$ is a local martingale. Also, 
\[\Var'[(W_S) \mid \mathcal{F}_{R_S}] ~ = ~ \Var'[(X) \mid \mathcal{F}_{R_S}] \cdot D_S ~ = ~ a[(X)]_{R_S} \cdot\frac{1}{a[(X)]_{R_S}}
~ = ~ 1 \enspace .
\]
So $(W)$ is a Brownian motion relative to $(\mathcal{F}_R)$, by Theorem~\ref{theorem:levy}.

Viewing the continuous processes $(X)$ and $(W)$ as measurable functions from $\Omega$ to $C_0([0,\infty))$,
the definition  $(W) := (X_R)$ can be written as $(W) = \Phi_a \circ (X)$, with probability law $P_{(W)} = (\Phi_a)_*(P_{(X)})$.
Since $(W)$ is a Brownian motion, $P_{(W)} = \mu_W$. So $\Phi_a$ is a measure-preserving bijection from
$(C_0([0,\infty)), P_{(X)})$ to $(C_0([0,\infty)), \mu_W)$. Hence its inverse $\Phi_a^{-1}$ is also measure-preserving, meaning
$P_{(X)} = (\Phi_a^{-1})_*(\mu_W)$.
\end{proof}

Theorems~\ref{theorem:stochastic-integral} and~\ref{theorem:unique} guarantee the existence and uniqueness up to equality-in-law of solutions $(X_t)$ to stochastic differential equations of the form
\[
X_0  =  0 \qquad \Exp'[(X) \mid \mathcal{F}]  =  0 \qquad \Var'[(X) \mid \mathcal{F}]  =  a[(X)] \enspace,
\]
for suitable adaptations $a$. 
One would of course like to have existence and uniqueness results more general 
stochastic differential equations of the form
\[
X_0  =  C_0 \qquad \Exp'[(X) \mid \mathcal{F}]  =  b[(X)] \qquad \Var'[(X) \mid \mathcal{F}]  = a[(X)] \enspace,
\]
where $C_0$ is an $\mathcal{F}_0$-measurable random variable specifying the initial distribution, and $a,b$ are 
suitable \emph{$C([0,\infty))$-adaptations} given as  functions
$C([0,\infty)) \times [0,\infty) \to \Real$. 
Because of Proposition~\ref{proposition:local-martingale-problem}, such questions amount to the \emph{well-posedness} of a 
\emph{local martingale problem}, so one can import known results from the literature, such as those 
of~\cite{SV}. In fact, our time-change 
proof of Theorem~\ref{theorem:unique} was influenced by the proof of \cite[Theorem~6.5.4]{SV}, although the setting is different.
(In \emph{op.\ cit.}, only Markov processes are considered, for which $a$ and $b$ depend just on the current $X_t$ rather than on the trajectory of $X$ over $[0,t]$, but the results are established for multi-dimensional Markov processes valued in $\mathbb{R}^d$, whereas, in the present section, we have only considered one-dimensional processes.)

%

\section{Discussion}
\label{section:end-notes}

In this section, we set  the present work in context with respect to other derivative-based approaches to stochastic calculus, beginning with a discussion of definitions based on right derivatives.
One such approach is taken by Nelson in \cite[Chapter 11]{Nelson}, where he  
directly adopts~\eqref{equation:drift-as-rd} as his definition of drift, 
which he calls the \emph{mean forward derivative}.
Nelson's analogue of variance rate is defined similarly, although in his case it is defined using \eqref{equation:nelson} (at deterministic times)  in place of  conditional variance  \eqref{rd:var}. To 
address the well-definedness issue raised in Section~\ref{section:intro}, Nelson assumes that $X_s$ is integrable at 
every time $s$. This allows him to use $\mathcal{L}^1$-convergence as his choice of limit. The third issue raised in Section~\ref{section:intro}, that of obtaining a useful stochastic calculus,  is not addressed. For example, the setting is too restrictive to derive  It\^o's formula. Nelson does, however, obtain a fundamental theorem of calculus, similar to our Theorem~\ref{theorem:stochastic-integral}, but valid in a much more restrictive context.


In later literature, e.g.~\cite[\S5.1]{RWII}, in which the rigorous development of stochastic calculus is based entirely on the stochastic integral,  one nonetheless finds 
formulas describing right derivatives for conditional expectation~\eqref{rd:exp} and for the projected-centre version of conditional variance~\eqref{equation:projected-centre}
offered to provide derivative-based intuitions for the concepts of drift and variance rate, without them being intended to be taken seriously as mathematical definitions. 

The historic literature on stochastic calculus, prior to It\^o's work, contains alternative derivative-based proposals for defining drift-like and variance-rate-like quantities, in which cut-offs are imposed to avoid integrability assumptions. For example, Kolmogorov~\cite{Kolm} and Doeblin~\cite{Doeblin} essentially define a version of drift as
\[
\lim_{t \downarrow s} \,\frac{\Exp[  \One_{(-1,1)} \cdot  (X_t -X_s) \mid \mathcal{F}_s ]}{t-s}  \enspace,
\]
although they do this only in the context of Markov processes, in which  conditioning on $ \mathcal{F}_s$ amounts to conditioning on $X_s$.
Doeblin, in particular, succeeds in using this definition, together with an analogous definition of a variance-rate-like quantity, to derive an analogue of
the It\^o formula, very similar to our Theorem~\ref{theorem:ito-formula}, that works at the level of generality of Markov processes satisfying the additional regularity conditions that Doeblin imposes. 

In the present article, the  simple device of defining derivatives in terms of limiting properties with respect to stopping times avoids integrability assumptions
without any need for explicit cut-offs, and without imposing additional regularity conditions on processes.
The resulting framework of stopping derivatives seems well adapted to establishing the core properties of stochastic calculus. It also provides a derivative-based perspective on several prominent concepts from stochastic calculus, such as local martingales (Theorems~\ref{theorem:martingale-driftless} and~\ref{theorem:driftless-martingale}), local martingale problems 
(Proposition~\ref{proposition:local-martingale-problem}) and  the It\^o integral (Theorem~\ref{theorem:stochastic-integral}).  
Theorem~\ref{theorem:stochastic-integral} shows that all processes obtained as stochastic integrals (of right-continuous processes with respect to Brownian motion) 
possess drift and variance rate. Thus our definitions of drift and variance rate are applicable to  It\^o processes
in essentially their full generality. At the same time, a potential advantage of our stopping-derivative-based definitions of drift and variance rate is that they can be applied locally to arbitrary processes at individual finite stopping times,  irrespective of whether or not  processes are given globally as It\^o integrals.

It is possible that stopping derivatives may elucidate other notions from the literature on stochastic processes.
For example, in the theory of diffusions, the \emph{infinitesimal generator} of a diffusion $(X)$ is the operator that 
maps a function $f$ to 
\begin{equation}
\label{equation:infinitesimal}
x \mapsto \lim_{t \downarrow 0} \frac{\Exp^x [f(X_t) - x ]}{t} \enspace 
\end{equation}
when defined, where $\Exp^x$ is the expectation for $(X)$ as a Markov process started at $X_0 = x$.
Dynkin's \emph{characteristic operator} is defined similarly, but with the limit over time $t$ replaced by a limit over 
\emph{first-exit times}, namely over stopping times
\[
T_{\varepsilon,x}(\omega) ~ : = ~ \inf \{ t \mid |X_t(\omega)- x | \geq \varepsilon\}  \enspace ,
\]
determined by $\varepsilon > 0$. The resulting characteristic operator maps a function $f$ to
\begin{equation}
\label{equation:characteristic} 
x \mapsto 
\begin{cases}
\lim_{\varepsilon \downarrow 0} \frac{\Exp^x [f(X_{T_{\varepsilon,x}}) - x]}{\Exp^x[T_{\varepsilon,x}]} & \text{if $\Exp^x[T_{\varepsilon.x}] < \infty$ for some $\epsilon > 0$} \\
0 & \text{otherwise} \enspace ,
\end{cases}
\end{equation}
again if this is defined.  The characteristic operator agrees with the  infinitesimal generator whenever the latter is defined, but is applicable in general to a broader class of functions. 
It can therefore be viewed as a 
generalisation of the 
infinitesimal generator obtained by switching from a time-based form of derivative~\eqref{equation:infinitesimal} to a space-based notion~\eqref{equation:characteristic}.
In the light of the present work, it seems plausible that the  characteristic operator can be given 
an alternative definition, using a stopping limit,  as the operator that maps a function $f$ to 
\begin{equation}
\label{equation:new-characteristic}
x \mapsto  \lim_{T \downarrow 0} \frac{\Exp^x [f(X_{T}) - x]}{\Exp^x[T]} \enspace.
\end{equation}
This proposed redefinition of the characteristic operator  
restores time as the axis of reference, with the resulting transition from~\eqref{equation:infinitesimal} to~\eqref{equation:new-characteristic} 
amounting to nothing more than a randomisation of the time variable.

An important characteristic of our definitions of drift and covariance rate as stopping limits is that they
take place locally in space and (random) time. Accordingly, it seems likely that these 
notions should generalise straightforwardly to stochastic processes on manifolds. 
This could lead to simple characterisations of Brownian motions on manifolds as processes
with zero drift and constant variance rate (cf.\ Theorem~\ref{theorem:levy}) and of local martingales 
on manifolds as processes with zero drift (cf.\ Theorem~\ref{theorem:driftless-martingale}). More generally, the locality of stopping derivatives may  offer other benefits as a potential formalism for expressing
stochastic differential equations on manifolds in a straightforward differential form. 

For processes valued in $\mathbb{R}^n$, 
it would be interesting to investigate whether the generalisation of It\^o's formula to path-dependent \emph{non-anticipative functionals} in~\cite[Theorem 4.1]{CF} can be recast in terms of stopping derivatives. 
Another natural question is whether stopping derivatives can be adapted to  processes, such as fractional Brownian motion, that do not have quadratic variation, but which fall into the scope of rough path theory~\cite{Lyo}.
Finally, we reiterate the point mentioned in Section~\ref{section:intro}
that the use of  derivatives with respect stopping times
has allowed us to bypass the issue 
that  $\Exp[X_t \mid \mathcal{F}_s]$ and
$\Var[X_t \mid \mathcal{F}_s]$ may diverge at deterministic times $t$. It is natural to wonder whether 
an analogous randomisation of time  might help to address other divergence problems that arise in mathematics and physics. 

\section*{Acknowledgements}

I thank Mihael Perman and Matija Vidmar for discussions and for helpful comments that led to mathematical improvements to the material presented in this article.

 \end{document}